\documentclass[12pt,twoside,reqno]{amsproc}
\usepackage[centertags]{amsmath}
\usepackage{amssymb,amsfonts,amsthm}
\newcommand{\TryPackage}[3]{\IfFileExists{#1.sty}{\usepackage{#1} #2}{#3}}
\TryPackage{mathrsfs}{\renewcommand{\mathcal}{\mathscr}}{%
   \TryPackage{eucal}{}{}}

\usepackage{a4ML}             % a4 paper format
\usepackage{math40,local}   % local makro definitions
\renewcommand{\aendrg}{\relax}

\renewcommand{\subsectionmark}[1]{\relax}

\newcommand{\signature}{\relax} % will be changed if doc class article is used
\begin{document}
\title[Boundary value problems I. Regularity and self--adjointness]{On boundary value problems for Dirac type operators. I. Regularity
  and self--adjointness}
% author one information
\author{Jochen Br\"uning}
\address{Humboldt Universit\"at zu Berlin\\
Institut f\"ur Mathematik\\ 
Unter den Linden 6\\
D--10099 Berlin
}
%\curraddr{}
\email{bruening@mathematik.hu-berlin.de}
\urladdr{http://spectrum.mathematik.hu-berlin.de}

%\thanks{}

% author two information
\author{Matthias Lesch}
\address{Humboldt Universit\"at zu Berlin\\
Institut f\"ur Mathematik\\ 
Unter den Linden 6\\
D--10099 Berlin
}
\curraddr{Universit\"at Bonn\\ Mathematisches Institut\\ Beringstr. 1\\
D--53115 Bonn}
\email{lesch@mathematik.hu-berlin.de, lesch@math.uni-bonn.de}
\urladdr{http://spectrum.mathematik.hu-berlin.de/$\sim$lesch}
\thanks{Both authors were supported by the Sonderforschungsbereich 288
of Deutsche Forschungsgemeinschaft.
The second named author was supported by the Gerhard Hess Programm and
a Heisenberg fellowship of Deutsche Forschungsgemeinschaft.}
%    General info
%\subjclass{Primary 58G20 boundary value problems on manifolds;
%Secondary 35F15 boundary value problems for linear first order PDE}

\date{May 28, 1999}

\begin{abstract} In a series of papers, we will develop systematically
the basic spectral theory of (self--adjoint) boundary value problems
for operators of Dirac type. We begin in this paper with the characterization
of (self--adjoint) boundary conditions with optimal regularity, for
which we will derive the heat asymptotics and index theorems
in subsequent publications. Along with a number of new results,
we extend and simplify the proofs of many known theorems.

Our point of departure is the simple structure which is displayed
by Dirac type operators near the boundary. Thus our proofs
are given in an abstract functional analytic setting, generalizing
considerably the framework of compact manifolds with boundary.

The results of this paper have been
announced in \cite{BruLes:STBVPDTO}.
\end{abstract}

\maketitle
\tableofcontents
\clearpage

\renewcommand{\thesubsection}{\thesection.\Alph{subsection}}
% well, we like experiments........
\section{Introduction}
\subsection{The compact case}\label{SecOneA}
Let $\tilde{M}$ be a compact oriented Riemannian manifold of dimension $m$, 
$\tilde{E} \to \tilde{M}$ a hermitian vector bundle, 
and $\tilde{D}$ a symmetric elliptic differential operator of order 
$d\in\Z_+$ on $C^\infty (\tilde{E})$. The classical results on existence, 
uniqueness,  and regularity of solutions of $\tilde{D}$ are among the 
cornerstones of Global Analysis:
\begin{itemize}
\begin{MLnumlist}
\MLitem{$ \tilde{D}$ is essentially self-adjoint in $ L^2(\tilde{E})$ 
with domain $C^\infty(\tilde{E})$ 
(by slight abuse of notation, we denote the closure of 
$\tilde{D}$ by the same symbol);}\eqnnumlbl{intro-1.1}
\MLitem{
$\tilde{D}$ is a Fredholm operator (of index $0$),
i.e. there are a bounded operator $\tilde{Q}$ and compact operators 
$\tilde{K}_r , \tilde{K}_l$ in $L^2(\tilde{E})$ such that}\eqnnumlbl{intro-1.2}
\[
 \tilde{D}\tilde{Q} = I - \tilde{K}_r,\qquad \tilde{Q}\tilde{D} = I-\tilde{K}_l ;
\]
\MLitem{with respect to the Sobolev scale 
$H_s(\tilde{E}):= {\cd}((\tilde{D}^2 + I)^{s/2d}),$ $s\ge 0,$ $\tilde{Q}$
is of order $-d$  and $\tilde{K}_{r/l} $ of order $- \infty.$}
\eqnnumlbl{intro-1.3}
\end{MLnumlist}
\end{itemize}

The restriction to symmetric operators is not essential since we may 
always consider a given elliptic operator together with its adjoint. 
But it is a technical advantage for more refined questions like 
index theorems: We bring in an isometric involution, $\tilde\omega,$ 
on $\tilde E$ which anticommutes with $\tilde{D}$ on $C^\infty (\tilde{E})$ 
and hence produces a splitting
\begin{equation}
 \tilde{D} = \begin{pmatrix} 0 & \tilde D_- \\ \tilde D_+ & 0\end{pmatrix}
 \text{ on } \; C^\infty (\tilde{E}_+) \oplus  C^\infty (\tilde{E}_-),
  \label{intro-1.4}
\end{equation}
with $\tilde{E}_\pm$ the $\pm 1$--eigenbundle of $\tilde{\omega}$. 
Then, by the well--known formula of McKean-Singer we have
\begin{equation}
 \ind \tilde{D}_+ = \tr \bigl[ \tilde{\omega} e^{-t\tilde{D}^2}\bigr],
  \quad t>0.\label{intro-1.5}
\end{equation}
It is hence of great importance that 
(cf. eg. \cite[Lemma 1.9.1]{Gil:ITHEASIT2E})
\begin{itemize}
\begin{MLnumlist}
\MLitem{
for any differential operator, $\tilde{P}$, of order $p$ on 
$C^\infty(\tilde{E})
$, we have an asymptotic expansion}
\begin{equation}
\tr \Bigl[ \tilde{P} e^{-t\tilde{D}^2}\Bigr] \sim_{t\to 0+} 
  \sum a_j (\tilde{D},\tilde{P}) t^{(j-m-p)/d}.
  \label{intro-1.6}
\end{equation}
\end{MLnumlist}
\end{itemize}

Even though it took a long time after the original proof of the 
Atiyah--Singer Index Theorem \cite{AtiSin:IEOCM}, \cite{Pal:SASIT}
until a complete proof could be based on \myref{intro-1.4} 
and \myref{intro-1.5} (cf. \cite{AtiBotPat:HEIT}, \cite{Gil:CELEC}), 
the heat equation method now seems to be the most powerful tool for 
extensions of the Index Theorem. Recall that for the important class 
of twisted Dirac operators, with $\tilde{E} = \tilde{S}\oplus\tilde{F}, 
\tilde{S}$ a spin bundle on $\tilde{M}$, this theorem reads
\begin{itemize}
\item $\DST\ind \tilde{D}_+ = \int\limits_{\tilde{M}} \hat{A} (\tilde{M}) 
\wedge \rm ch \,\tilde{F}.$\eqnnumlbl{intro-1.7}
\end{itemize}

\subsection{Compact manifolds with boundary}\label{SecOneB}
In this series of papers, we want to present the extension of the main 
results quoted above (\myref{intro-1.1}, \myref{intro-1.2}, 
\myref{intro-1.3}, \myref{intro-1.6}, \myref{intro-1.7}) 
to Dirac type operators on manifolds with boundary.

Though a good part 
of our results is more or less known, we obtain a conceptually as well as 
technically transparent derivation of this theory, with considerable 
simplifications and extensions in most cases. Moreover, the functional 
analytic approach we have developed lends 
itself naturally to substantial generalizations, e.g. to situations 
with non-compact boundaries. The basic inspiration for this approach
is, of course, the beautiful work of Atiyah, Patodi, and Singer 
\cite{AtiPatSin:SARG} which we generalize.

To explain our work in greater detail, we consider a compact hypersurface,
$N$, in $\tilde{M}$  which bounds an open subset, $M$, of $\tilde{M}$. 
We assume that $N$ is oriented as the boundary of $M$.
We put 
$E:=\tilde{E} \restr M,$ \linebreak[3]
$D:=\tilde{D} \restr C^\infty (\overline{M})$. \linebreak[3]
Then $D$ is a first order elliptic differential operator on $M$, 
and symmetric in $L^2(E)$ with domain $C^\infty_0(E)$. 
Recall that $D$ is an operator of \emph{Dirac type}
\cite[Def. 3.36]{BerGetVer:HKDO},
\cite[Sec. 1.8.2]{Gil:ITHEASIT2E}
if $D^2$ has scalar principal symbol given by the
metric tensor, i.e.
$\hat D(\xi)^2=\lvert\xi\rvert^2$ for each $\xi\in T^*M$.
Note that this class of operators is considerably larger
than the class of \emph{Dirac operators associated to a Clifford
connection} (or twisted Dirac operators) 
\cite[Sec. II.5]{LawMic:SG}, \cite[p. 119]{BerGetVer:HKDO}.

If $D$ is of Dirac type
then we obtain in a tubular neighbourhood, $U$, of $N$ 
in $\tilde{M}$ a very simple separation of variables. In fact, $U$ is 
isometric to $(-\varepsilon_0, \varepsilon_0) \times N$ with metric 
$dx^2\oplus g_N(x), x\in (-\varepsilon_0, \varepsilon_0),$
$g_N$ a smooth family of metrics on $N$, and with 
$E_N := \tilde{E} \restr N$ we obtain the following result.
\begin{lemma}\label{intro-S1} 
Let $D$ be of Dirac type. As operator in $L^2(E\restr U)$ with domain 
$C^\infty_0(E\restr U)$, $D$ is unitarily equivalent to an operator of the form
\begin{equation}
\gamma\left(\frac{d}{dx} +A(x)\right)+V(x)
\label{intro-1.8}
\end{equation}
in $L^2((-\varepsilon_0,\varepsilon_0),L^2(E_N))$
with domain $C^\infty_0((-\varepsilon_0,\varepsilon_0),
C^\infty(E_N)).$ 

Here, $\gamma\in{\cl}(L^2(E_N))$ and $A(x)$ is (the closure of) 
a symmetric elliptic differential operator on $E_N$ of first order; 
${\cd}(A(x))=: {\cd}$ is independent of $x$ and $A(x)$ depends smoothly 
on $x\in (-\varepsilon_0,\varepsilon_0).$ Furthermore,
$V\in\cinf{(-\eps_0,\eps_0),\cinf{\operatorname{End}(E_N)}}.$

Moreover, the following relations hold:
\begin{subequations}\label{intro-1.9}
\begin{align}
  &\gamma^\ast =-\gamma,\quad \gamma^2=-I,\label{intro-1.9a}\\
  &\gamma({\cd})= {\cd} \text{ and } \gamma A(x) + A(x)\gamma = 0, 
\quad x\in(-\varepsilon_0,\varepsilon_0).\label{intro-1.9b}
\end{align}
\end{subequations}
If $D$ is a Dirac operator associated to a Clifford connection, then
$A$ can be chosen in such a way that $V=0$.
\end{lemma}
This lemma has been widely used for some time, especially in the product case 
$(g_N(x)\equiv g_N(0))$ where it plays a prominent role in \cite{AtiPatSin:SARG}. 
For non-product metrics some care is needed to compute $A(x)$ in each specific
case,
cf. e.g. \cite[Sec. 3.10]{Gil:ITHEASIT2E} and 
\cite[Sec. 5]{Bru:LITCMROT}.

We will base our analysis on a thorough study of the operator equation 
\myref{intro-1.8} with the structure properties \myref{intro-1.9};
\emph{these properties will be assumed throughout this paper}.
This approach 
is reasonable since the results we are aiming at can be obtained from 
merging "interior analysis" (to be carried out on $\tilde{M}$) with 
"boundary" analysis involving the operator \myref{intro-1.8}.

The main difference between the analysis of $D$ and $\tilde{D}$ lies, of
course, in the fact that $D$ is not  essentially self-adjoint on
$C^\infty_0(E).$ 
Moreover, if self-adjoint extensions of $D$ exist, they may differ widely with 
respect to existence, uniqueness, regularity, and heat trace expansions. 
It is, therefore, our first task to characterize those self-adjoint extensions 
which behave nicely with respect to existence, uniqueness, and regularity; 
this is the purpose of the present paper.

\subsection{Results for the model operator}\label{SecOneC}
Replacing in \myref{intro-1.8} $L^2(E_N)$ by an arbitrary Hilbert space, $H$, 
and $C^\infty(E_N)$ by the domain, $H_1$, of  a self-adjoint operator $A$ in
$H$, we obtain the model operator
\begin{equation}
 D=\gamma\left(\frac{d}{dx} +A\right) \text{ in } \,{\ch}_0:=L^2(\R_+, H)
   \text{ with domain } C^\infty_0(\R_+, H_1).\label{intro-1.10}
\end{equation}

We will have to deal with variable coefficients but for the purpose of the
present introduction we will restrict to the constant coefficient case.
Indeed, for most of the problems dealt with in this paper operators with
variable coefficients merely appear as perturbations of \myref{intro-1.10},
in view of the Kato-Rellich Theorem. Furthermore, since $V(x)$ is 
a bounded operator, it can be ignored in the discussion of self--adjoint
extensions of $D$. 

On $C^\infty_0(\R_+, H_1)$ we clearly have
\begin{equation}
   (Df, g)_{{\ch}_0} - (f, Dg)_{{\ch}_0} = \scalar{f(0)}{\gamma g(0)}_H.
   \label{intro-1.11}
\end{equation}
Now if $D$ is symmetric on a subspace, ${\cd}^0$, of $C^\infty_0(\R_+, H_1)$
then it follows from \myref{intro-1.11} that, with $I-P$ the
orthogonal projection onto $\ovl{\cd^0}$ in $H$, we have
\begin{subequations}\label{intro-1.12}
\begin{equation}
 {\cd}^0 \subset {\cd}_P := \bigsetdef{f\in C^\infty_0 (\R_+, H_1)}{P f(0)=0}
  \label{intro-1.12a}
\end{equation}
and
\begin{equation}
  I - P \le \gamma^\ast P\gamma.\label{intro-1.12b}
\end{equation}\end{subequations}
Moreover, 
$D_{P,0} := D \restr {\cd}_P$ is a symmetric extension of 
$D\restr \cd^0$. If we assume for a moment that $H$ is of
finite dimension then it is readily seen that $D_{P,0}$ 
is essentially self-adjoint in ${\ch}_0$ if and only if
\begin{equation}
  I - P = \gamma^\ast P\gamma.\label{intro-1.13}
\end{equation}
Indeed, if $D_{\max}$ denotes $(D\restr C^\infty_0 ((0,\infty), H_1))^\ast$ then
\begin{equation}
  {\cd}(D_{\max}) \subset H_{1,\loc} (\R_+, H),\label{intro-1.14}
\end{equation}
and \myref{intro-1.11} remains valid with $D_{\max}$ in place of $D$, 
for $f, g \in {\cd}(D_{\max}).$

An orthogonal projection, $P$, with \myref{intro-1.13} will be called 
$\gamma$--\emph{symmetric}; it is easy to see that $\gamma$--symmetric 
projections -- and hence self--adjoint extensions of $D$ -- exist 
if and only if 
\begin{equation}
  \sign (i \gamma \restr \ker A) = 0.\label{intro-1.15}
\end{equation}
As an illustration, note that $i\frac{d}{dx}$ does not admit self-adjoint 
extensions in $L^2(\R_+).$

Returning to the general case, we meet the essential difficulty that 
\myref{intro-1.14} has no reasonable analogue. In particular, elements of 
${\cd}(D_{\max})$ do not  admit $H$-valued restrictions to zero. To overcome
this obstacle, we imitate the Sobolev scales $H_s(E_N)$ and $H_s(E)$ and their
interplay in our abstract setting (which has some tradition in
Analysis, cf. eg. \cite[Chap. XIII]{Pal:SASIT}). $H_s(E_N)$ is replaced by
\begin{equation}\begin{split}
 H_s&:= H_s(A)\\
    & := \left\{\begin{array}{l}
  {\cd}(|A|^s), \text{ equipped with the graph norms for } s\ge 0; \\
\text{a suitable dual of } H_{-s} (A), \text{ for } s<0.\end{array}\right.
		\end{split}
\end{equation}
We also need
\begin{align*}
  H_\infty&:= H_\infty(A):=\bigcap_{s\in\R} H_s(A);\\
\intertext{and}
  H_{-\infty}&:=H_{-\infty}(A):=\bigcup_{s\in\R} H_s(A).
\end{align*}
Next we introduce, for $n\in\Z_+$,
\begin{subequations}\label{intro-1.17}
\begin{equation}
 {\ch}_n := {\ch}_n (\R_+, A) := \bigcap^n_{k=0} H_k (\R_+ , H_{n-k}(A)),
 \label{intro-1.17a}
\end{equation}
where, for $i,j\in\Z_+ ,$
\begin{equation}\begin{split}
  H_i(&\R_+, H_j (A)) \\
     &:= \bigsetdef{f\in{\ch}_0}%
{\Bigl(\frac{d}{dx}\Bigr)^l f\in L^2(\R_+, H_j(A)), \quad 0\le l\le i}.
  \label{intro-1.17b}
		\end{split}
\end{equation}
\end{subequations}
By interpolation, we then obtain a scale of Hilbert spaces, 
${\ch}_s ={\ch}_s (\R_+, A), s\in \R_+.$ \myref{intro-1.17a},
\myref{intro-1.17b} also make sense with $\R$ in place of $\R_+$;
in this way we obtain the scale $\ch_s(\R,A)$. 
Generalizing the classical Trace Theorem for Sobolev spaces, 
we have the following result about trace maps which will allow 
the formulation of boundary conditions.
\begin{theorem}\label{intro-S1a}
The map 
$$
 r : C^\infty_0 (\R_+, H_\infty)\ni f \longmapsto f(0) \in H_\infty 
$$
extends by continuity to a map
$$
  r_s : {\ch}_s \longrightarrow H_{s-1/2}, \quad s> 1/2,
$$
and also to a map
$$ 
 r^\ast : {\cd} (D_{\max}) \to H_{-1/2} .
$$
\end{theorem}
Of course, the loss of regularity under the trace map requires the 
(continuous) extension of the boundary projections to the space $H_{-1/2}$. 
To deal with this, we introduce \emph{operators of finite order} on the 
Hilbert scale $(H_s(A))_{s\in\R}$. 
Thus, a linear map, $B:H_\infty\to H_\infty,$ is an operator of order
$\mu\in\R$ if for each $s\in\R$ there is a constant $C(s)$ such that, for any
$x\in H_\infty$,
\begin{equation}
 \| Bx\|_{H_s} \le C(s) \| x\|_{H_{s+\mu}} \,.\label{intro-1.18}
\end{equation}
In particular, $B$ extends to an element of ${\cl}(H_s, H_{s-\mu})$ 
for all $s\in\R.$ The totality of such operators forms the linear space 
$\Op^\mu(A).$ $\Op^{-\infty}(A) := \bigcap_{\mu\in\R}\Op^\mu (A)$ 
is called the space of \emph{smoothing operators}.

Thus we will have to require that the boundary projections are elements of 
$\Op^0(A)$. It follows easily from \myref{intro-1.18} that $\Op^0(A)$ is a 
$\ast$-algebra but it is, in general, not \emph{spectrally invariant} 
in the sense that $B\in \Op^0(A)$ and $B$ invertible in ${\cl}(H)$ 
implies $B^{-1}\in \Op^0(A)$.
% (cf. Proposition \plref{I.S1-2.5} below). 
To allow for a minimum of functional constructions, we do need actually 
even more. We are forced to restrict attention to certain subalgebras, 
$\Psi^0(A) \subset \Op^0(A),$ satisfying the following two conditions. 

%\begin{itemize}
%\item[($\Psi$1)] $\Psi^0(A)$ is a $\ast$-subalgebra of
%$\OpAn$ with \emph{holomorphic functional calculus};\\
%\item[($\Psi$2)] $\Psi^0(A)$ contains an orthogonal projection,
%$P_+(A)$, satisfying 
%\begin{equation}
% I-P_+(A)=\gamma^*P_+(A)
%\gamma,\quad P_{(0,\infty)} (A)\le P_+(A)\le P_{[0,\infty)}(A).
% \label{intro-1.19}
%\end{equation}
%\end{itemize}
\begin{alphnumbering}[axiomsintro]
\begin{MLlist}
\begin{MLnumlist}
\MLitem{$\PsiA$ is a $\ast$-subalgebra of
$\OpAn$ containing the smoothing operators and
with holomorphic functional calculus;} \puteqnnum\label{intro-1.19b}

\medskip
\MLitem{$\PsiA$ contains an orthogonal projection
$P_+(A)$, satisfying}
\begin{equation}
 I-P_+(A)=\gamma^*P_+(A)
\gamma,\quad P_{(0,\infty)} (A)\le P_+(A)\le P_{[0,\infty)}(A).
 \label{intro-1.19}
\end{equation}
\end{MLnumlist}
\end{MLlist}
\end{alphnumbering}

Recall that an algebra, ${\ca}\subset {\cl}(H),$ has holomorphic 
functional calculus if for $B\in {\ca}$ and $f$ holomorphic in a 
neighbourhood of $\spec B$ (in $\cl(H)$) we have $f(B) \in{\ca}$, 
where $f(B)$ is defined by the Cauchy integral. 
Note also that the existence of $P_+(A)$ with \myref{intro-1.19} 
is equivalent to 
\myref{intro-1.15} in the finite dimensional case.
In general, if $0 \notin \specess A$ then \myref{intro-1.15} is equivalent to 
the existence of a spectral projection of $A$ satisfying \myref{intro-1.19}.

Let us illustrate these conditions for the case where $A$ is an elliptic
differential 
operator on $C^\infty(E_N)$ and $N=\partial M$ as in Sec. \ref{SecOneB}.
Then we have $H_s(A) \simeq H_s(E_N)$, and a 
natural choice of the algebra $\Psi^0(A)$ is the algebra of classical 
pseudodifferential operators on $E_N$, to be denoted by 
$\Psi^0_{\class} (E_N)$. It follows from results of Seeley 
\cite[Thm. 5]{See:CPEO} that $\Psi^0_{\class}(E_N)$ has 
holomorphic functional 
calculus. Moreover, since $0 \notin \specess A$, we have to verify 
\myref{intro-1.15} to obtain a spectral projection, $P_+(A)$, of $A$ 
fulfilling \myref{intro-1.19}; but this is a consequence of the 
Cobordism Theorem.
To see this, 
we split $H=: H_+\oplus H_-$ according to the $\pm i$--eigenspaces of
$\gamma$. In view of 
\myref{intro-1.9},
\begin{equation}
  A=\begin{pmatrix} 0 & A_- \\ A_+ & 0 \end{pmatrix},
   \label{I.G1-1.17}
\end{equation}
and $A_+$ is a Fredholm operator with index
\begin{equation}
 \ind A_+ = \dim \ker A\cap \ker(\gamma-i) -\dim \ker A \cap \ker (\gamma + i).
 \label{I.G1-1.18}
\end{equation}
Now it is straightforward to check that there exists 
an orthogonal projection $P_+(A) \in \OpAn$ with the property
\myref{intro-1.19} if and only if
\begin{equation}
  \ind A_+ =0,
  \label{I.G1-1.19}
\end{equation}
and this follows from the Cobordism Theorem
(cf. the discussion after \cite[Cor. 3.6]{BruLes:EICNLBVP}). 

Again from Seeley's work, we deduce that 
$P_+(A) \in\Psi_{\class}^0(E_N)$ so \myref{intro-1.19b}, \myref{intro-1.19}
are satisfied 
in a natural way.

Now we are in the position to formulate our results for the model operator. 
The main theorem of this article reads as follows.
\begin{theorem}\label{intro-S2}
Let $H$ be a Hilbert space and $A$ a self-adjoint operator in $H$. 
Assume, moreover, that an algebra, $\Psi^0(A)\subset \Op^0(A),$ 
is given with the properties \myref{intro-1.19b} and 
\myref{intro-1.19}.

Then $D_{P,0}$ with domain \myref{intro-1.12a} is essentially self-adjoint in 
$L^2(\R_+, H)$ for any orthogonal projection $P\in \Psi^0(A)$ with the 
properties
\begin{equation}
 \gamma^\ast P\gamma = I - P\tag{\ref{intro-1.13}}
\end{equation}
and
\begin{equation}
  (P,P_+(A)) \text{ is a Fredholm pair.}\label{intro-1.20}
\end{equation}
The domain of the closure, $D_P$, of $D_{P,0}$ is 
\begin{equation}
   \cd(D_P)=\bigsetdef{f\in\ch_1(\R_+,A)}{Pf(0)=0}.
   \label{I.G7-1.24}
\end{equation}

Conversely, if $A$ is discrete then the self--adjointness of 
$D^*\restr\cd(D_P)$ %$\bigsetdef{f\in\ch_1(\R_+,A)}{Pf(0)=0}$
implies \myref{intro-1.13} and \myref{intro-1.20}.
\end{theorem}

We note that the orthogonal projection $P_-(A)=I-P_+(A)\in\Psi^0(A)$
obviously does not satisfy \myref{intro-1.20}. However, we
will show in Proposition \plref{I.S1-4.8} 
that $D_{P_-(A)}$ is essentially self--adjoint with domain
$\cd(D_{P_-(A)})\supsetneqq\bigsetdef{f\in\ch_1(\R_+,A)}{Pf(0)=0}$
(cf. Proposition \plref{I.S1-4.7} through \plref{I.S1-5.10} 
for a detailed discussion of this phenomenon).
Hence the "self--adjointness" in the last statement of the Theorem
cannot be replaced by "essentially self-adjoint on $\cd_P$".

The crucial notion of a Fredholm pair of projections is described in 
Sec. \ref{Secthree}; 
the proof of Theorem \plref{intro-S2} takes up Sections \ref{Sectwo} to 
\ref{Secfive}; 
since we do not see a direct way to prove it, we have to interpolate 
various notions of "regularity" which accounts for the length of our 
presentation. At the end of Section \ref{Secfive} we give the proof of 
Theorem \plref{intro-S2} refering to the several intermediate
results.

We can view Theorem \plref{intro-S2} as the analogue of \myref{intro-1.1}
for the model operator. Taking advantage of the self-adjointness of $D_P$ 
we can try to satisfy \myref{intro-1.2} by setting
\begin{equation}
  Q := \int\limits_{|\lambda|\ge 1} \lambda^{-1} d  E(\lambda),
  \label{intro-1.21}
\end{equation}
where $E(\lambda) = E_{D_P}(\lambda), \lambda\in \R,$ denotes the spectral 
resolution of $D_P$. 

The regularity result in Theorem \plref{I.S5-4.13}
together with the compactness property expressed in Proposition 
\plref{I.S1-2.13} easily yields the following analogue of 
\myref{intro-1.2} and \myref{intro-1.3}.
\begin{theorem}\label{intro-S3}
We assume the situation of Theorem \plref{intro-S2} and, in addition, 
that $A$ is discrete. For $\phi\in C^\infty_0(\R)$ with $\phi=1$ near 
$0$ we put $Q_\phi := \phi  Q$. Then $Q_\phi$ maps into ${\cd}(D_P)$ 
and there are compact operators, $K_{r/l,\phi}$, 
in ${\ch}$ such that
\begin {equation}
  D_PQ_\phi = \phi - K_{r,\phi},\qquad 
  Q_\phi D_P = \phi -K_{l,\phi}.
  \label{intro2-1.22}
\end{equation}
Moreover, $Q_\phi $ is of order $-1$ and $K_{r/l,\phi}$ of order $-\infty$ 
with respect to the Sobolev scale ${\ch}_s (\R_+, A), s\in\R_+.$
\end{theorem}

\subsection{Results for manifolds with boundary}\label{SecOneD}
It is not difficult to translate Theorems 
\plref{intro-S2} and \plref{intro-S3}
into statements on $D$ and the Sobolev scales $H_s(E)$ and 
$H_s(E_N), s\in \R$. We only need to make Lemma \plref{intro-S1} somewhat 
more explicit.

For this, we introduce, on $U$, the global coordinate
\begin{equation}
 x(p) := \operatorname{dist}(p,N) , \quad p\in U,
 \label{intro-1.23}
\end{equation}
and denote by
\begin{equation}
 \Phi : L^2(E\restr U) \longrightarrow 
  L^2 ((-\varepsilon_0, \varepsilon_0), L^2(E_N))
  \label{intro-1.24}
\end{equation}
the isometry implicit in Lemma \plref{intro-S1}. Then we have the properties
\begin{subequations}\label{intro-1.25}
\begin{align}
&\psi\Phi u= \Phi \bigl((\psi\circ x)u\bigr),
   &&\psi \in C^\infty (-\varepsilon_0, \varepsilon_0), u\in L^2(E\restr U),
              \label{intro-1.25a}\\ 
    & (\Phi u)(0) = u \restr N, & &u\in C_0^\infty (E\restr U),
               \label{intro-1.25b}\\
   &  \Phi\bigl((\psi\circ x) H_s(E)\bigr) =  
      \psi{\ch}_s(\R,A), 
   &&\psi\in C^\infty_0 (-\varepsilon_0,\varepsilon_0), s\in \R,
\end{align}
\end{subequations}
which allow us to localize near $N$ and to transfer regularity.

To formulate the boundary conditions, we restrict attention to orthogonal 
projections in $L^2(E_N)$ which are classical pseudodifferential operators 
i.e. from now on we choose 
$\Psi^0(A) =\Psi^0_{\class} (E_N),$ as indicated above. 
In the theory of boundary value problems for linear elliptic differential 
operators, it was observed by Calder\'on \cite{Cal:BVPEE} 
that a prominent role is played by an idempotent, 
$C^+ \in\Psi^0_{\class} (E_N)$, with the property that
\begin{equation}\begin{split}
 C&^+(H_s(E_N)) = N_s(E_N)\\
  & := \bigsetdef{u\in H_s(E_N)}%
{u=\tilde{u}\restr N \text{ for } \tilde{u}\in H_{s+1/2}(E)
\text{ with } D\tilde{u} = 0}\label{intro-1.26}
		\end{split}
\end{equation}
for all $s\in\R; C^+ $ is called the \emph{Calder\'{o}n projector} 
(cf. \cite{See:SIBVP,See:TPO}, a comprehensive summary can be found in
\cite[Appendix]{Gru:TEPBPDTOMGS}). 
One checks that
\begin{equation}
 C^+ - P_+(A) \in \Psi^{-1}_{\class}(E_N)
\end{equation}
which explains the importance of the Atiyah--Patodi--Singer boundary condition. 
In order to obtain boundary conditions which define Fredholm operators 
(as in \myref{intro-1.2}, \myref{intro-1.3}), Seeley introduced the notion 
of "well-posed" boundary condition \cite{See:TPO}
which we develop in Section \ref{Secseven}\ below. 
Combining the results described so far with Theorem \plref{I.S1-7.2}, 
we obtain the following optimal version of Seeley's result, as a consequence 
of our general theory.
\begin{theorem}\label{intro-S4} Let $M$ be a compact manifold with
boundary as in Section \plref{SecOneB} 
and let $D$ be an operator of Dirac type (resp. a
first order symmetric elliptic differential operator of the
form \myref{intro-1.8}) acting on sections of the hermitian vector
bundle $E$.\aendrg

Let $P\in\Psi^0_{\class} (E_N)$ be an orthogonal projection in 
$L^2(E_N)$ satisfying \myref{intro-1.13}.

Then $D_P:=D^*\restr {\cd}_P 
:= \bigsetdef{f\in H_1(\tilde{E}\restr\overline{M})}%
{P(f\restr N) = 0}$ is self-adjoint in $L^2(E)$ if and only if 
$P$ is well-posed in the sense of Seeley. This, in turn, is equivalent
to the fact that $(P,P_+(A))$ is a Fredholm pair.

\aendrg
In this case there are a bounded operator, $Q$, and compact operators, 
$K_{r}, K_{l},$ in $L^2(E)$ such that 
$Q$ maps into ${\cd}(D_P)$ and 
$$
  D_PQ = I-K_r, \quad QD_P= I-K_l.
$$
Moreover, with respect to the Sobolev scale $H_s(E), s\in \R_+, Q$ 
is of order $-1$ and $K_{r/l}$ of order $-\infty$.
\end{theorem}

The proof of this Theorem is presented at the end of Section
\plref{Secseven}.

\subsection{Further results}\label{SecOneE}
Theorem \plref{intro-S4} is the main application presented here of the 
results in this paper but, by their abstract character, they apply to more 
singular situations as well. This will be carried out in part III of 
this series for covering spaces of compact manifolds with boundary.

Among the material presented here we still have to mention 
Sec. \ref{Secsix}\ where 
we deal with variable coefficients and also, in preparation for the following 
parts, with the regularity theory of $D^2$. Here the reader will find the 
(fairly easy) arguments necessary to prove the results mentioned above for 
variable coefficients (i.e. for the case where $g_N(x)$ is not constant 
near $x=0$).

The other publications in this series will be devoted to the analogues of 
the statements \myref{intro-1.6} and \myref{intro-1.7}.

In part II, we will develop systematically the index theory of the 
operator\linebreak[3] 
$\gamma\bigl(\frac{d}{dx} + A(x)\bigr)$ (with suitable involutions) on 
finite or infinite intervals, with appropriate boundary conditions. 
Again, we will simplify and extend various known results but also present 
some new theorems.

Part III will address index theorems on manifolds with boundary
in full generality, based 
on Theorem \plref{intro-S4}, and we will apply it to the 
"glueing" of indices. Moreover, we use our techniques to give very simple 
proofs of various results in index theory like the Cobordism Theorem, 
or the index theorems of Callias and Ramachandran.

Part IV will derive the heat expansion for the operators described in 
Theorem \plref{intro-S4}, using only the simple structure of the model 
operator. In this context, variable coefficients present
essential new difficulties. 
We will pay special attention to the (spectrally defined) 
determinant of these operators as a function on the Grassmannian of 
well-posed projections.

Some of the results of this and the forthcoming papers have been
announced in \cite{BruLes:STBVPDTO}.

%  Projekt einer Reihe von (bis zu) vier Arbeiten mit JB
%  ueber Randwertprobleme fuer Dirac Operatoren
%
%  Teil I. Regularitaetstheorie

%  Der eigentliche TeX Text:

%  Letzte Aenderung: 28.05.99

\section[Sobolev spaces and operator algebras]{Sobolev spaces and
operator algebras associated with self--adjoint operators}
\label{Sectwo}

\subsection{The Sobolev scale of an unbounded operator}

We consider a Hilbert space, $H$. We fix a self--adjoint
(unbounded) operator $A$ in $H$ and introduce the dense subspace
\begin{equation}
  H_\infty:=\bigcap_{s\ge 0} \cd(|A|^s),
 \label{I.G1-2.1} 
\end{equation}
where $\cd$ denotes the domain of an unbounded operator. For $s\in\R$
let $H_s=H_s(A)$ be the completion of $H_\infty$ with respect to the scalar
product
\begin{equation}
   \scalar{x}{y}_s:= \scalar{(I+A^2)^{s/2}x}{(I+A^2)^{s/2}y},
   \label{I.G1-2.2}
\end{equation}
such that $H_0=H$, $H_1=\cd(A)$, and we have embeddings
$i_{s',s}:H_{s'}\subset H_s$ of norm at most $1$, for $s'>s$.
Then $(I+A^2)^{\mu/2}$ induces an isometry $H_s\to H_{s-\mu}$ for
all $s,\mu\in \R$ and we obtain a perfect pairing
\begin{equation}\begin{split}
   &H_{-s}\times H_s\longrightarrow \C,\\
   &(x,y)\longmapsto \Bs{x}{y}:=
    \scalar{(I+A^2)^{-s/2}x}{(I+A^2)^{s/2}y}_0,
  		\end{split}
   \label{I.G1-2.3}
\end{equation}
with $\Bs{x}{y}=\scalar{x}{y}_0$ for $x,y\in H_\infty$.
In particular
\begin{equation}
     |\Bs{x}{y}|\le \|x\|_{-s}\|y\|_s.
    \label{I.G2-2.4}
\end{equation}
If $s\ge 0$ then 
\begin{equation}
   H_s=\cd(|A|^s).
   \label{I.G1-2.4}
\end{equation} 
 
The natural inclusion map
$H_{s'}\hookrightarrow H_s$ is compact 
if and only if $A$ is \emph{discrete}, i.e. if $A$ has a compact
resolvent in $H$.
For $s\in\R$, $H_s$ is a Banach space with strong antidual $H_{-s}$.
$H_\infty$ becomes a Fr\'echet space under the 
seminorms $(\|\cdot\|_n)_{n\in\Z}$, with dual space
$H_{-\infty}=\cup_{s\in\R} H_s$. 

\begin{dfn}\label{I.S1-2.1} A linear map, $T:H_\infty\to H_\infty$,
is called an \emph{operator  of order $\mu\in\R$}, if $T$ 
induces a continuous linear map
from $H_s$ to $H_{s-\mu}$, for all $s\in \R$, that is, if for any
$s\in\R$ there is a constant $C_s(T)$ such that
\begin{equation}
     \|Tx\|_{s-\mu}\le C_s(T) \|x\|_s,\quad x\in H_\infty.
   \label{I.G1-2.5}
\end{equation}
$T$ is called \emph{smoothing} or \emph{of order} $-\infty$ if it
is of order $\mu$ for all $\mu\in \R$.
We denote the set of all operators of order $\mu$ by $\OpA{\mu}$,
$-\infty\le \mu<\infty$.
\end{dfn}
\begin{prop}\label{I.S1-2.2}
\begin{thmenum}
\item For $T\in \OpA{\mu}$ let $T^*$ be the Hilbert space adjoint
of $T$ considered as an (unbounded) operator in $H$ with domain 
$H_\infty$. 
Then $H_\infty\subset\cd(T^*)$ and 
$T^t:=T^\ast\restr H_\infty$ is in $\OpA{\mu}$.
In particular,
$T^*$ is densely defined (resp. bounded if $\mu\le 0$).
\item Suppose that $T,T^t:H_\infty\to H_\infty$ are linear with 
  \[ \scalar{T x}{y} = \scalar{x}{T^t y}, \quad x,y\in H_\infty .\]
If $T$ and $T^t$ satisfy \myref{I.G1-2.5} for some $s_j\ge 0, j\in\Z_+$, with
$\lim\limits_{j\to\infty} s_j=\infty$, 
then both $T$ and $T^t$ are in $\OpA{\mu}$. 
\item Let $T\in\cl(H_0)$  with $T(H_s)\subset H_{s-\mu}, T^*(H_s)\subset 
H_{s-\mu}$
for $s\ge 0$ and some fixed $\mu\le 0$. Then $T\in\OpA{\mu}$.
\end{thmenum}
\end{prop}
\begin{proof} (1)
We note first that \myref{I.G1-2.5} is equivalent to the estimate
\begin{custeqnnum}{\ref{I.G1-2.5}'}
\begin{equation}
|\scalar{Tx}{y}| \le C_s(T) \|x\|_s \; \|y\|_{\mu-s} ,\quad x,y\in H_\infty, s\in \R.
   \label{I.G2-2.5'}
\end{equation}
\end{custeqnnum}
This implies (1), with $C_s(T^\ast) = C_{\mu-s}(T).$

(2) In (2), we have \myref{I.G2-2.5'} for $s_j, j\in \Z_+$, 
by assumption on $T$, and from the assumption on $T^t$ 
we derive \myref{I.G2-2.5'} for $-s_j+\mu, j\in \Z_+$. 
Hence the assertion follows from complex interpolation.

(3) It follows from the Closed Graph Theorem that
$T$ and $T^*$ map $H_s$ continuously into $H_{s-\mu}$, for $s\ge 0$.
Then (2) implies the assertion. 
\end{proof}

Note that for any rapidly decreasing Borel function $f:\R\to \R$ the
operator $f(A)$ is smoothing. In particular,
$P_0(A)$, the orthogonal projection onto $\ker A$, is smoothing.

\begin{cor}\label{I.S1-2.3} $\OpAn$ is a $*$--subalgebra of $\cl(H_0)$,
the algebra of bounded linear operators on $H_0$.

The smoothing operators $\OpA{-\infty}$ form a two--sided
$*$--ideal in $\OpAn$.
\end{cor}
\begin{proof} The product of operators of order $0$ is an operator
of order $0$. Thus $\OpAn$  is an algebra and by the previous
Proposition it is a $*$--algebra.
The last statement is obvious.\end{proof}

\begin{lemma}\label{I.S1-2.4}
Let $T\in \Op^\mu (A).$ 
If $T(H_\infty)$ is finite-dimensional, then $T$ is smoothing.
%$T\in \OpA{-\infty}.$
\end{lemma}
\begin{proof} 
Choose a basis $(f_i)^N_{i=1} $ for $T(H_\infty)$ which is orthonormal in $H$. 
We can write, for $x\in H_\infty$,
  \[ Tx = \sum^N_{i=1} T_i (x) f_i \,,\]
with $T_i: H_\infty\to \C$ linear. It follows for $s\in \R$ that 
  \[|T_i(x)| = |\scalar{Tx}{f_i}| \le C_s (T) \|x\|_s \; \|f_i\|_{\mu-s} \,. \]
But then $T_i(x) = \scalar{x}{e_i}$ for some $e_i\in H_\infty.$
\end{proof}

\begin{dfn}\label{I.S3-2.5}
Slightly more generally, we call a family
$(H_s)_{s\in\R}$ $(s\in\R_+)$ a \emph{scale of Hilbert spaces} if
\begin{thmenum}
\item \label{Hscaleone}
$H_s$ is a Hilbert space for each $s\in\R$ ($s\in\R_+$),
\item\label{Hscaletwo}
$H_{s'}\hookrightarrow H_s$ embeds continuously for $s\le s'$,
\item\label{Hscalethree}
if $s < t, 0<\theta <1, $ then the complex interpolation space satisfies
\[[H_s,H_t]_\theta =H_{\theta t+(1-\theta)s},\]
\item\label{Hscalefour}
$H_\infty:=\bigcap\limits_{s\in\R_+} H_s$ is dense in $H_t$ for each $t$.
\end{thmenum}

In view of \myref{I.G1-2.3} the Sobolev scale of an 
unbounded self--adjoint operator satisfies in addition
\begin{thmenum}\addtocounter{enumi}{4}
\item \label{Hscalefive}
the $H_0$--scalar product restricted to $H_\infty$ extends to
an antidual pairing between $H_s$ and $H_{-s}$ for all $s\in\R$.
\end{thmenum}
\end{dfn}

For the complex interpolation
method we refer to \cite[Sec. 4.2]{Tay:PDEI}.
A scale $(H_s)_{s\in\R_+}$ satisfying \myref{Hscaleone}--\myref{Hscalefour}
can be extended
to a scale of Hilbert spaces parametrized over $\R$
by defining $H_{-s}$ to be the completion
of $H_\infty$ with respect to the norm
\[ \|x\|_{-s}:=\sup_{y\in H_\infty\setminus\{0\}}
    \frac{\lvert \scalar{x}{y}\rvert}{\|y\|_{s}}, \quad s\in\R_+.\]
This family will also satisfy \myref{Hscalefive}. However, the scales
of Sobolev spaces on
manifolds with boundary usually do not
satisfy \myref{Hscalefive}.

Thus, if $A$ is a self-adjoint operator in the Hilbert space $H$, then 
$(H_s(A))_{s\in\R}$ is a scale of Hilbert spaces satisfying
\myref{Hscalefive}.
The converse is almost true: namely, if $(H_s)_{s\in\R}$ is a scale
of Hilbert spaces satisfying \myref{Hscalefive},
then for each $N>0$ there exists a self--adjoint
operator $A$ in $H_0$ with $H_s(A)=H_s$ for $|s|\le N$ (cf.
\cite[Sec. 1.2.1]{LioMag:NHBVPAI}). However, it is not clear
whether there exists such an $A$ for all $s\in\R$ simultaneously.
For example, we will prove in Corollary \plref{I.S2-2.16}
that $(\ch_s(\R_+,A))_{s\in\R_+}$ satisfies 
\myref{Hscaleone}--\myref{Hscalefour}, hence it can be extended to
a scale of Hilbert spaces parametrized over $\R$. 
However, we do not know of a self--adjoint
operator $B$ in $L^2(\R_+,H)$ such that $\ch_s(\R_+,A)=H_s(B)$ for all
$s\ge 0$. 

The results of this section have obvious counterparts for a
given scale of Hilbert spaces $(H_s)_{s\in\R}.$ 
If the axiom \myref{Hscalefive} is not satisfied, duality
arguments are not possible and one has to restrict attention
to $s\ge 0$.

A useful property of an algebra 
${\ca}\subset {\cl}(H)$ is its \emph{spectral invariance} 
by which we mean the assertion
\begin{equation} 
T\in{\ca}, T \; \text{invertible in} \,{\cl}(H) \Longrightarrow T\; 
\text{invertible in }{\ca} . 
\end{equation}
A slightly stronger assumption is
\begin{equation}
  {\ca} \;\text{admits holomorphic functional calculus,}
\end{equation}
by which we mean that for $T\in {\ca}$ and 
any function $f$, holomorphic in a neighborhood of $\spec T$, 
we have $f(T) \in {\ca}.$

We proceed to show that $\Op^0(A)$ is, in general, not spectrally
invariant; 
this will force us to restrict our attention to suitable subalgebras.

\begin{prop}\label{I.S1-2.5} Let $A$ be a 
discrete operator with eigenvalues
$|\mu_0|\le |\mu_1|\le |\mu_2|\le\ldots \to\infty$. Assume that there
exists a subsequence $(\mu_{n_k})_{k\in\Z_+}$ satisfying
\begin{equation}
   0<\gl_1\le \Bigl|\frac{\mu_{n_k}}{\mu_{n_{k+1}}}\Bigr|\le \gl_2<1
   \label{I.G1-2.8}
\end{equation}
for some $\gl_1,\gl_2$.

Then $\OpAn$ is not spectrally invariant in $\cl(H_0)$, i.e.
there exists an operator $T\in\OpAn$ which is invertible in $H_0$
but with $T^{-1}\not\in\OpAn$.
\end{prop}

\begin{proof} Denote by $(e_n)_{n\in\Z_+}$ an orthonormal basis of
$H$ with $A e_n=\mu_n e_n$. Put for $\xi=\sum\limits_{n=0}^\infty \xi_n e_n\in
H$
\[ K\xi:= \sum_{k=0}^\infty \xi_{n_{k+1}} e_{n_k}.\]
$K$ is in $\OpAn$ in view of \myref{I.G1-2.8}. Furthermore,
we have
\[ \|K\|_{H_0\to H_0}=1.\]
Hence for any $0<\gl<1$ the operator $I+\gl K$ is invertible
in $H_0$.

Next consider $\xi=(\xi_n)_{n\in\Z_+}$ with
\[\xi_n:=\begin{cases}
    (-1)^k \gl^{-k}, &n=n_k,\\
    0                & \text{otherwise.}
    	 \end{cases}\]
\myref{I.G1-2.8} implies $\xi\in H_s$ for $s\le s_0<0$ small enough.
Since $(I+\gl K)\xi=0$ the operator $(I+\gl K)$ is not invertible in
$H_s$
and thus $(I+\gl K)^{-1}\not\in \OpAn.$\end{proof}

\begin{remark}\label{I.S1-2.6} 
\begin{enumerate}
\item Gramsch \cite[Ex. 6.2]{Gra:RISOA} gave
the first example of a discrete operator $A$ such that
$\OpAn$ is not spectrally invariant. Proposition \plref{I.S1-2.5} is a
generalization of this.
\item We conjecture that for discrete
$A$ the algebra $\OpAn$ is never spectrally invariant in $\cl(H_0)$.
\item The condition \myref{I.G1-2.8} is fulfilled if
\[ \mu_n\sim C n^\ga\]
for some $C,\ga>0$. One just takes $n_k:=2^k$. Thus, if
$A$ is a self--adjoint elliptic operator on a compact manifold
then $\OpAn$ is not spectrally invariant.
\end{enumerate}
\end{remark}

\subsection{Sobolev spaces on $\R$ and $\R_+$}

Let 
\begin{equation}
 D=\gamma(\frac{d}{dx}+A)
 \label{I.G1-2.9}
\end{equation}
as defined in \myref{intro-1.9}, \myref{intro-1.10}.
%, with domain
%$\cinfz{(0,\infty),H_\infty}$ in $L^2(\R_+,H)$.

Since it will be necessary to distinguish
between operators on the whole line $\R$ and on the half line
$\R_+$ we denote by
$\Dtilde$ the operator $\gamma(\frac{d}{dx}+A)$ acting on
$\cinfz{\R,H_\infty}$ in $L^2(\R,H)$. 
This is a symmetric operator in $L^2(\R,H)$ 
with $\tilde D^2=-\frac{d^2}{dx^2}+A^2$.
Furthermore, we put $D_\pm:=\tilde D\restr\cinfz{\R_\pm^*,H_\infty},
\R_\pm^*:=\R_\pm\setminus\{0\}$. $D_\pm$ is a symmetric operator
in $L^2(\R_\pm,H)$.
If no confusion is possible
we will write $D:=D_+$. 

The unique self--adjoint extension of $\Dtilde$ will be the source
of another Sobolev scale which is crucial for our further study.

\begin{lemma}\label{I.S1-2.8} All powers of $\Dtilde$ are essentially self--adjoint.
\end{lemma}
\begin{proof}  The assertion is easy to see for bounded $A$.
But $\Dtilde$ commutes with the spectral projections of $|A|$ so we can
reduce the problem to this case.\end{proof}

By a slight abuse of notation, we identify $\Dtilde$ with
its unique self--adjoint extension.
Next we define the Sobolev spaces
\begin{equation}
  \ch_s(\R,A):=H_s(\Dtilde),\quad s\in\R.
  \label{I.G1-2.10}
\end{equation}
$\ch_s(\R,A)$ is a Hilbert space with norm
\begin{equation}
\|f\|_s^2=\int_\R \|(I+\xi^2+A^2)^{s/2} \hat f(\xi)\|^2 \dbar \xi,
 \label{I.G1-2.11}
\end{equation}
where
\[\hat f(\xi):=\int_\R 
    e^{-ix\xi} f(x) dx,\quad
  \dbar \xi:=\frac{1}{2\pi} d\xi.
\]
The following maps are continuous:
\begin{equation}\begin{split}
    \Dtilde, \frac{d}{dx}, A&:\ch_s(\R,A)\longrightarrow \ch_{s-1}(\R,A),\\
    \gamma, \check{\hphantom{a}}&:\ch_s(\R,A)\longrightarrow
     \ch_s(\R,A),
   		\end{split}
    \label{I.G2-2.12}
\end{equation}
where $(Af)(x):=A(f(x)), (\gamma f)(x):=\gamma f(x), \check{f}(x):=f(-x).$

Furthermore, we put $\ch_s(\R_\pm,A):=\{f\restr \R_\pm\,|\, f\in
  \ch_s(\R,A)\}$, $s\in\R$, and equip this space with the quotient
Hilbert space structure, i.e.
\begin{equation}
   \|f\|_s:=\inf\bigsetdef{\|\tilde f\|_s}{\tilde f\restr \R_\pm=f}.
  \label{I.G1-2.12}
\end{equation}
Clearly, there is a natural restriction map $\ch_s(\R,A)\to
\ch_s(\R_\pm,A)$ of norm $1$. \myref{I.G2-2.12} also holds with $\R_\pm$
in place of $\R$, except that $\check{\hphantom{a}}$ maps $\ch_s(\R_\pm,A)$
continuously into $\ch_s(\R_\mp,A)$. 

We will prove in Corollary \plref{I.S2-2.16}
below that the family $(\ch_s(\R_+,A))_{s\in\R_+}$ is a scale
of Hilbert spaces hence, as noted 
after Definition \plref{I.S3-2.5}, it can be extended to a scale
of Hilbert spaces parametrized over $\R$. 
We do not claim, however, that $(\ch_s(\R_+,A))_{s\in\R}$ is
a scale of Hilbert spaces. Although this is conceivable, the family
$(\ch_s(\R_+,A))_{s\in\R}$ will certainly not satisfy axiom
\myref{Hscalefive}. Hence, in the discussion of $\ch_s(\R_+,A)$ we will
restrict ourselves to the case $s\ge 0$. According to the remark
after Definition \plref{I.S3-2.5} we could extend $(\ch_s(\R_+,A))_{s\ge 0}$
to a scale of Hilbert spaces satisfying \myref{Hscalefive}. However,
we refrain from doing so for two reasons: first, this definition
of Sobolev spaces of negative order would differ from the usual
definition on manifolds with boundary \cite[Sec. 4.5]{Tay:PDEI} and,
secondly, integration by parts shows that the operator $D$ would
not be of order one for that scale. The spaces
$\ch_s(\R_+,A)$ for negative $s$ will play no further role in the rest of
the paper.

We note that for $s\ge 0$ we have an inclusion
$\ch_s(\R_{(+)},A)\subset L^2(\R_{(+)},H_s)$ of norm less or equal 1.
Indeed, 
\begin{align}
    \int_\R \|f(x)\|_{H_s}^2 dx&=\int_\R \|\hat f(\xi)\|_{H_s}^2 \dbar\xi=
       \int_\R \|(I+A^2)^{s/2} \hat f(\xi)\|^2\dbar\xi
       \notag\\
     &\le \int_\R \|(I+\xi^2+A^2)^{s/2}\hat f(\xi)\|^2\dbar\xi
            \label{I.G1-2.13}\\
     &= \|f\|_s^2.\notag
\end{align}

\begin{remark}\label{I.S1-2.9}
The Sobolev spaces introduced above can be expressed in terms
of standard Sobolev spaces. Namely, we put for a Hilbert space $H$
\begin{equation}\begin{split}
  H_k&(\R_{(+)},H)\\
     &:=\bigsetdef{f\in L^2(\R_{(+)},H)}{
    \bigl(\frac{d}{dx}\bigr)^j f\in L^2(\R_{(+)},H), 0\le j\le k}\\
    &=H_k(\R_{(+)})\hat\otimes H,
		\end{split}
   \label{I.G2-2.14}
\end{equation}
where $\hat\otimes$ denotes the Hilbert space tensor product.
Then one infers from \myref{I.G2-2.12} that for
$n\ge 0$ 
\begin{equation}
    \ch_n(\R,A)=\bigcap_{k=0}^n H_k(\R,H_{n-k}(A)).
     \label{I.G1-2.14}
\end{equation}
(cf. \cite[p. 8]{RobSal:SFMI}).
\myref{I.G1-2.14} is also true with $\R_+$ in place of $\R$ 
(this is shown in Corollary \plref{I.S2-2.16} below), 
this fact, however, is less obvious.

\myref{I.G1-2.14}
can be improved using complex interpolation. 
We consider Hilbert spaces
$E,F$ with self--adjoint operators $B\ge I$ in $E$ and $C\ge I$ in $F$.
We put $E':=\cd(B)$, $F':=\cd(C)$. Then $B\otimes C$ is a symmetric
operator on $H_\infty(B)\otimes H_\infty(C)\subset E\hat\otimes F$.
This operator is essentially self--adjoint and its unique self--adjoint
extension will be denoted by $B\hat\otimes C$ (cf.
\cite{BruLes:RTPUOHC}). Since $B\ge I$ we can define the graph
norm of $B$ by $\|x\|_B:=\|Bx\|_F$. If we define the graph norms
of $C$ and $B\hat\otimes C$ similarly we see that 
\begin{equation}
      \cd(B\hat\otimes C)= E'\hat\otimes F'.
     \label{I.G3-2.18}
\end{equation}
Furthermore, $B\hat\otimes I$ and $I\hat\otimes C$ are commuting
self--adjoint operators with $(B\hat\otimes I)(I\hat\otimes C)=B\hat\otimes C$.
Hence we have for any $\theta>0$
\begin{equation}
     (B\hat\otimes C)^\theta= B^\theta\hat\otimes C^\theta.
     \label{I.G3-2.19}
\end{equation}
These  considerations imply
\begin{equation}\begin{split}
       [E\hat\otimes F,E'\hat\otimes F']_\theta&=
            \cd((B\hat\otimes C)^\theta)=\cd(B^\theta\hat\otimes C^\theta)\\
             &= \cd(B^\theta)\hat\otimes \cd(C^\theta)
      =[E,E']_\theta \hat\otimes [F,F']_\theta.
   \end{split}   \label{I.G3-2.20}
\end{equation}

From the inequality
$\frac 12 (b^s+c^s)\le (b+c)^s\le 2^s(b^s+c^s)$,  for $b,c,s\ge 0$, and
the Spectral Theorem we infer
\begin{equation}\begin{split}
    [E,E']_s\hat\otimes F\cap E\hat\otimes [F,F']_s&=
        \cd(B^s\hat\otimes I)\cap\cd(I\hat\otimes C^s)\\
    &= \cd((B\hat\otimes I+I\hat\otimes C)^s)\\
    &= [E\hat\otimes F,E'\hat\otimes F\cap E\hat\otimes F']_s.
  \end{split}
      \label{I.G3-2.21}
\end{equation}

With these preparations we can prove:

\begin{prop}\label{I.S5-2.10}
For $s\ge 0$ we have the identities
\begin{align}
 \bigcap_{0\le t\le s} H_t(\R_{(+)},H_{s-t}) &= 
    H_s(\R_{(+)}, H_0)\cap H_0(\R_{(+)},H_s),
   \label{I.G2-2.21}\\
\intertext{and}
{\ch}_s (\R,A) &= H_s(\R,H_0)\cap H_0(\R,H_s).\label{I.G2-2.22}
\end{align}
\myref{I.G2-2.22} is to be understood as an equality of Hilbert spaces,
i.e. we have norm estimates
\begin{equation}\begin{split}
C^{-1}\bigl(&\norm{f}_{H_s(\R,H_0)}+\norm{f}_{H_0(\R,H_s)}\bigr)\\
            &\le \norm{f}_{\ch_s(\R,A)}\le 
C\bigl(\norm{f}_{H_s(\R,H_0)}+\norm{f}_{H_0(\R,H_s)}\bigr)
\label{I.G6-2.24}
		\end{split}
\end{equation}
for all $f\in\chs$.
\end{prop}
\begin{proof}
We apply \myref{I.G3-2.18}--\myref{I.G3-2.21} to the spaces
$H_s(\R_{(+)},H_t)=H_s(\R_{(+)})\hat\otimes H_t$.
A consequence of \myref{I.G3-2.20} is the identity
\begin{equation}
[H_s(\R_{(+)} , H_t), H_{s'}(\R_{(+)} , H_{t'})]_\theta = H_{\theta
s'+(1-\theta)s}(\R_{(+)}, H_{\theta t'+(1-\theta)t}) ,
    \label{I.G2-2.19}
\end{equation}
in particular
\begin{equation}
[L^2(\R_{(+)} , H_s), H_s(\R_{(+)},H_0]_\theta = 
   H_{\theta s} (\R_{(+)} , H_{(1-\theta) s}).
   \label{I.G2-2.20}
\end{equation}
This implies \myref{I.G2-2.21}. 
In view of \myref{I.G1-2.14}, \myref{I.G3-2.21},
\myref{I.G2-2.20}, and \myref{I.G2-2.21}
we obtain for  $0\le s\le N\in\Z_+$ 
\begin{align}
{\ch}_s (\R,A) &= [L^2(\R,H_0), {\ch}_N (\R,A)]_{s/N}\notag\\
         &=  [L^2(\R,H_0), H_N(\R,H_0)\cap H_0(\R,H_N)]_{s/N}\notag\\
 &= [L^2(\R,H_0), H_N(\R,H_0)]_{s/N}\cap[L^2(\R,H_0),H_0(\R,H_N)]_{s/N}\notag\\
 &= H_s(\R,H_0)\cap H_0(\R,H_s).
\end{align}
The first inequality in \myref{I.G6-2.24} follows from \myref{I.G1-2.13}
and a similar calculation for $\norm{\cdot}_{H_s(\R,H_0)}$. Since
both norms $\norm{\cdot}_{H_s(\R,H_0)}$, 
$\norm{\cdot}_{H_s(\R,H_0)}+\norm{\cdot}_{H_0(\R,H_s)}$ are Hilbert
space norms on $\ch_s(\R,A)$, the first inequality in \myref{I.G6-2.24}
implies the second.
\end{proof}

The same reasoning would work for $\R_+ $ in place of $\R$ 
if we knew that $({\ch}_s(\R_+,A))_{s\ge 0}$ is a scale of Hilbert spaces; 
this fact is the content of Corollary \plref{I.S2-2.16}.
\end{remark}

\begin{lemma}\label{I.S1-2.7} Let $f:\R\times\spec A\to \C$ be a continuous
function. Assume that for fixed $x\in\R$ the function $f(x,\cdot)$ is
of polynomial growth. 
\begin{thmenum}
\item For $u\in H_\infty$, we have
\begin{equation}
\int_\R\|f(x,A)u\|^2dx\le \|u\|^2
    \sup_{\gl\in\spec A}\int_\R |f(x,\gl)|^2dx.
   \label{I.G3-2.26}
\end{equation}
\item For $\varphi\in \cinfz{\R,H_\infty}$,
the function $x\mapsto f(x,A)\varphi(x)$
is weakly integrable over $\R$ and
\begin{equation}
\| \int_\R f(x,A) \varphi(x) dx\|^2\le \|\varphi\|^2_{L^2(\R,H)}
     \sup_{\gl\in\spec A} \int_\R |f(x,\gl)|^2 dx.
   \label{I.G3-2.27}
\end{equation}
\end{thmenum}
By continuity, \textnormal{(1)} extends to $u\in H$ and \textnormal{(2)} extends
to $\varphi\in L^2(\R,H)$.
\end{lemma}
\begin{proof} $f(x,A)u$ is well--defined since $f(x,\cdot)$ is assumed to
be of polynomial growth.

(1)
Let $(E(\gl))_{\gl\in\R}$ be the spectral resolution of $A$.
Then
\begin{align*}
       \int_\R\|f(x,A)u\|^2dx&= \int_\R \scalar{f(x,A)u}{f(x,A)u}dx\\
     &= \int_\R\int_\R |f(x,\gl)|^2 d\scalar{E(\gl)u}{u} dx\\
   &\le  \|u\|^2\sup_{\gl\in\spec A}\int_\R |f(x,\gl)|^2dx.
\end{align*}

\smallskip
\noindent
(2) For $u\in H_\infty$ we estimate with (1):
\begin{align*}
   \int_\R |\scalar{f(x,A)\varphi(x)}{u}|dx &\le \int_\R \|\varphi(x)\|\,
        \|f(x,A)u\|dx\\
    &\le \|\varphi\|_{L^2(\R,H)}\biggl(\int_\R \|f(x,A)u\|^2 dx\biggr)^{1/2}\\
   &\le \|\varphi\|_{L^2(\R,H)}\, \|u\| 
    \biggl(\sup_{\gl\in\spec A} \int_\R |f(x,\gl)|^2 dx\biggr)^{1/2},
\end{align*}
which proves both the weak integrability and the estimate.\end{proof}

\begin{theorem}[Trace Theorem] \label{I.S1-2.10}
\begin{thmenum}
\item For $s>1/2$ 
the restriction 
map \[r:\cinfz{\R,H_\infty}\longrightarrow H_\infty, 
\quad f\longmapsto f(0),\] 
induces by continuity bounded linear operators
$r:\chs\to H_{s-1/2}$ and $r_+:\chsplus\to H_{s-1/2}$.
\item  Let $s>-1/2$. If
$f, (I+A^2)^{-1/2}f'\in\chs$, then we have the estimate
\[\|f(0)\|_{s-1/2}\le C (\|f\|_s+\|(I+A^2)^{-1/2} f'\|_s).\]
\end{thmenum}
\end{theorem}
\begin{proof} (1) For $f\in\cinfz{\R,H_\infty}$ we write $f(0)=\int_\R \hat
f(\xi)\dbar\xi$ and estimate with Lemma \plref{I.S1-2.7}:
\begin{align*}
    \|&f(0)\|_{s-1/2}^2=\Bigl\|\int_\R (I+A^2)^{s/2-1/4} \hat f(\xi)\dbar \xi\Bigr\|^2\\
    &= \Bigl\|\int_\R (I+A^2)^{s/2-1/4}(\xi^2+I+A^2)^{-s/2} 
      (\xi^2+I+A^2)^{s/2}\hat f(\xi)\dbar \xi\Bigr\|^2\\
  &\le \sup_{\gl>0} \int_\R (1+\gl)^{s-1/2}(1+\gl+\xi^2)^{-s}\dbar\xi\;
       \|f\|_s^2\\
  &=\int_\R (1+\xi^2)^{-s}\dbar\xi\; \|f\|_s^2.
\end{align*}
This proves the assertion for $\chs$.
For $\chsplus$, it is an immediate consequence of the definition of
the quotient norm on $\chsplus$.

(2) It suffices to prove this estimate for 
$f\in\cinfz{\R,H_\infty}$. By (1) we have
%\begin{epeqnarray}
\begin{align*}
   \|f&(0)\|_{s-1/2}^2= \|(I+A^2)^{-1/2}f(0)\|_{(s+1)-1/2}^2\\
        &\le C(s)\int_\R \|(I+A^2)^{-1/2}(I+\xi^2+A^2)^{1/2}
               (I+\xi^2+A^2)^{s/2}\hat f(\xi)\|^2d\xi\\
     &= C(s) \int_\R \|(I+\xi^2+A^2)^{s/2} \hat
        f(\xi)\|^2 d\xi \\
    &\quad + C(s) \int_\R \|(I+\xi^2+A^2)^{s/2}(I+A^2)^{-1/2}\xi 
        \hat f(\xi)\|^2 d\xi\\
    &= C(s)(\|(I+A^2)^{-1/2} f'\|_s^2+\|f\|_s^2).\tag*{\qed}
\end{align*}
\exendproof
%\end{epeqnarray} 

The restriction map $r$ is in fact surjective, more precisely:

\begin{prop}\label{I.S1-2.11} Let $u\in H_{s-1/2}$, $\varphi\in H_s(\R)$, and put
$T:=(I+A^2)^{1/2}$. Then
the function
\[f(x):=\varphi(xT)u\]
is in $\chs$ and
\[\|f\|_s\le \|\varphi\|_s \|u\|_{s-1/2}.\]

In particular, if $\varphi(0)=1$ we get a continuous right inverse
to the restriction  map $r$ by \[e_\varphi(u)(x):=\varphi(xT)u.\]
\end{prop}
Note that this result is valid without any restriction on $s$.

\begin{proof} Assume first $\varphi\in\cs(\R)$ and $u\in H_\infty$. 
Then  $f\in\cs(\R,H_\infty)$ and from the spectral theorem we
see that
\[\hat f(\xi)=T^{-1}\hat \varphi(\xi T^{-1})u.\]
With Lemma \plref{I.S1-2.7} (1) we obtain
\begin{align*}
  \|f\|_s^2&= \int_\R \| (\xi^2+T^2)^{s/2} T^{-1}\hat \varphi(\xi T^{-1})u\|^2
     \dbar\xi\\
   &\le \sup_{\gl>0} \int_\R (\xi^2+\gl^2)^s \gl^{-2s-1} |\hat \varphi(\xi/\gl)|^2
        \dbar\xi\; \|u\|_{s-1/2}^2\\
   &= \|\varphi\|_{H_s(\R)}^2 \|u\|_{s-1/2}^2.
\end{align*}
Since $\cs(\R)$ is dense in $H_s(\R)$ and $H_\infty$ in $H_{s-1/2}$,
we reach the conclusion.\end{proof}

\begin{cor}\label{I.S1-2.12} Let $s>1/2+k, k\in\Z_+$. Then the map
\[r_+^{(k)}:\chsplus\longrightarrow \moplus_{j=0}^{k}H_{s-1/2-j},\quad
   f\longmapsto \moplus_{j=0}^k (D^jf)(0),\]
is continuous and surjective. 
Moreover, there exists a continuous right inverse $e^{(k)}$ to
$r_+^{(k)}$ with $e^{(k)}\restr \moplus\limits_{j=0}^{k} H_\infty$
independent of $s$.
\end{cor}
\begin{proof} Clearly, by Theorem \plref{I.S1-2.10} $r_+^{(k)}$ exists
and is continuous since 
\[D^j:\ch_s(\R_+,A)\longrightarrow \ch_{s-j}(\R_+,A)\]
is. It remains to construct $e^{(k)}$.

We choose $\varphi \in
C^\infty_0(\R)$ with $\varphi=1$ near $0$. Then in view
of Proposition \plref{I.S1-2.11} we may choose $e^{(0)}:=e_\varphi$;
this is independent of $s>1/2$. 

Inductively, we assume that we have constructed $e^{(k)}$. Choose $\varphi\in
C^\infty_0(\R)$ with
\[
  \varphi^{(j)} (0) = 0, \quad 0\le j\le k,  \quad \varphi^{(k+1)} (0) =1,
\]
and put for $(\xi_0,\ldots ,\xi_{k+1}) \in \bigoplus\limits^{k+1}_{j=0}
H_{s-1/2-j}, s> 1/2+k+1,$
\begin{align*}
e^{(k+1)}&(\xi_0,\ldots ,\xi_{k+1})\\
&:= (-\gamma)^{k+1} \varphi(xT)T^{-k-1}
(\xi_{k+1}-(D^{k+1} e^{(k)}(\xi_0,\ldots ,\xi_k))(0))\\
&\quad+e^{(k)}(\xi_0,\ldots ,\xi_k),
\end{align*}
where $T=(I+A^2)^{1/2}$, as in Proposition \plref{I.S1-2.11}.
\end{proof}

\begin{lemma}\label{I.S1-4.1} 
\begin{thmenum}
\item $\cd(D^*)\subset \{f\in L^2(\R_+,H)\,|\, f'\in L^2(\R_+,H_{-1})\}$,
and we have a continuous restriction map
$r^*:\cd(D^*)\to H_{-1/2}, f\mapsto f(0)$.
\item For $f_j\in \cd(D^*), j=1,2$ we have
  \[-(D^*f_1,f_2)+(f_1,D^*f_2)=
  \lim_{\eps\to 0}\scalar{\gamma (\eps A+i)^{-1} f_1(0)}{% 
    (\eps A+i)^{-1} f_2(0)}.\]
If $f_j(0)\in H_{1/2}$ for some $j\in\{1,2\}$ then
\begin{equation}
 -(D^*f_1,f_2)+(f_1,D^*f_2)=\Bss{(-1)^{j+1}1/2}{\gamma f_1(0)}{f_2(0)}.
  \label{I.G1-4.1}
\end{equation}
\end{thmenum}
\end{lemma}
\begin{proof} (1) For $f\in\cd(D^*)$ we have
$Af\in L^2(\R_+,H_{-1})$, hence $f'\in L^2(\R_+, H_{-1})$.
We apply Theorem \ref{I.S1-2.10} (2) with $s=0$ and obtain the estimate
\begin{align*}
\|f(0)\|_{-1/2}&\le C (\|f\|+ \|(I+A^2)^{-1/2}f'\|)\\
      &\le C(\|f\|+\|(I+A^2)^{-1/2}D^*f\|+\|(I+A^2)^{-1/2}\gamma Af\|)\\
      &\le C(\|f\|+\|D^*f\|).
\end{align*}

\smallskip
\noindent
(2) For $f\in \cd(D^*)$ we put $f_\eps(x):=i(\eps A+i)^{-1} f(x)$.
Then we have $f_\eps\in\cd(D^*)\cap L^2(\R_+,H_1)$ and
$f_\eps'\in L^2(\R_+,H)$. Moreover, $f_\eps\to f, D^*f_\eps=(D^*f)_{-\eps}
\to D^*f$, as $\eps\to 0$.
Thus, integration by parts gives
\begin{align*}
   -(D^*f_1,f_2)+(f_1,D^*f_2)&= \lim_{\eps\to 0}\big[
   (D^*f_{1,\eps},f_{2,\eps})-(f_{1,\eps},D^*f_{2,\eps})\big]\\
    &=\lim_{\eps\to 0}\scalar{\gamma f_{1,\eps}(0)}{f_{2,\eps}(0)}\\
  &=\lim_{\eps \to 0}\scalar{\gamma (\eps A+i)^{-1} f_1(0)}{%
       (\eps A+i)^{-1} f_2(0)}.
\end{align*}

To prove the last assertion we note that
if $f_j(0)\in H_{1/2}$ then $\lim_{\eps\to 0}i(\eps A+i)^{-1}f_j(0)=f_j(0)$ in
the $H_{1/2}$--topology and we reach the conclusion.\end{proof}

A similar result holds for $D_-^*$. This follows immediately from
the identity
\begin{equation}
       D_-^*\gamma \check{f}=-\gamma (D_+^*f)^\vee.
      \label{I.G2-2.23}
\end{equation}
A consequence of the previous lemma is the following characterization
of the space $\ch_n(\R_+,A)$.

\begin{prop}\label{I.S2-2.13}  Let $f\in\cd((D_+^*)^n)$. Then $f\in\ch_n(\R_+,A)$
if and only if there exists a $g\in\cd((D_-^*)^n)$ with
\begin{equation}
(D_-^jg)(0)=(D_+^jf)(0),\quad  0\le j\le n-1.
   \label{I.G2-2.24}
\end{equation}
\end{prop}

\begin{proof} If $f\in\ch_n(\R_+,A)$ then there exists $\tilde f\in\ch_n(\R,A)$
with $\tilde f\restr\R_+=f$ and one can take $g:=\tilde f\restr \R_-$.

Conversely, assume that $g\in\cd((D_-^*)^n)$ satisfies \myref{I.G2-2.24}.
Put 
\begin{equation}
     \tilde f(x):=\begin{cases}
                   f(x),&  x\ge 0,\\
                   g(x),&  x<0.
            	  \end{cases}
     \label{I.G2-2.25}
\end{equation}
For $\varphi\in\cinfz{\R,H_\infty}$ we find using Lemma \plref{I.S1-4.1}
\[\begin{split}
  |\scalar{\tilde f}{\tilde D^n \varphi}|&
    =|\scalar{(D_+^*)^n f}{\varphi\restr \R_+}_{L^2(\R_+,H_0)}
       + \scalar{(D_-^*)^n g}{\varphi\restr \R_-}_{L^2(\R_-,H_0)}|\\[0.5em]
    &\le C_{\tilde f}\, \|\varphi\|_{L^2(\R,H_0)},
  \end{split}
\]
thus $\tilde f\in\cd((\Dtilde^n)^*)=\ch_n(\R,A)$ in view of Lemma
\plref{I.S1-2.8}.\end{proof}

From these facts we get the following regularity result.

\begin{cor}\label{I.S2-2.14} 
\textnormal{(1)} For $n\in\Z_+$ we have
\begin{equation}\begin{split}
  \ch_n&(\R_+,A)\\
     &=\biggl\{f:\R_+\to H\,\biggm |\, \begin{array}{l} 
       D^j f\in L^2(\R_+,H), \;0\le j\le n,\\
      (D^jf)(0)\in H_{n-1/2-j},\;0\le j\le n-1 
\end{array}\biggr\}.
		\end{split}
   \label{I.G2-2.18}
\end{equation}
\textnormal{(2)}
\[\ch_{n+1}(\R_+,A)=\bigsetdef{f\in L^2(\R_+,H)}{Df\in \ch_n(\R_+,A),\;
    f(0)\in H_{n+1/2}}.\]
\end{cor}
\begin{proof} (1) The inclusion "$\subset$" follows from \myref{I.G2-2.12} and
Corollary \plref{I.S1-2.12}. To prove the converse inclusion we pick
$f$ in the right hand side of \myref{I.G2-2.18} and put
\[g(x):= e^{(n-1)}(f(0), (Df)(0),\ldots, (D^{n-1}f)(0))(x), \quad x\le 0.\] 
Now the assertion follows from Proposition \plref{I.S2-2.13}.

(2) follows immediately from (1) by induction.
\end{proof}

\begin{prop}\label{I.S2-2.15}
For a fixed integer $N$ choose $(a_j)^N_{j=1} $ such that
\begin{equation}
\sum^N_{j=1} (-j)^l a_j = 1, \quad l=0,1,\ldots ,N-1.
   \label{I.G2-2.26}
\end{equation}
For $f\in C^\infty_0 ([0,\infty), H_\infty)$ put
\[
  (E^{(N)} f)(x) :=\begin{cases}
                   f(x),&  x\ge 0,\\
 \sum\limits^N_{j=1} a_j f(-jx),&x<0.
                   \end{cases}
\]
Then $E^{(N)}$ has a continuous extension to a map
\[ 
   H_s(\R_+,H_0)\cap H_0(\R_+,H_s)\longrightarrow
  {\ch}_s(\R,A), \quad\text{for}\quad  s\le N-1,
\]
also denoted by $E^{(N)}$.
\end{prop}
\begin{proof}
The system \myref{I.G2-2.26} determines the
$a_1, \ldots , a_N$ uniquely.

Fix $k\in\Z_+, k\le N-1$, and $f\in H_k(\R_+, H_0)\cap H_0(\R_+,H_k).$ 
Then \myref{I.G2-2.26} implies that the right and left
derivatives of $E^{(N)}f$ coincide up to order 
$k$ at $0$. Hence $E^{(N)}$ extends by continuity to a map
\[
   H_k(\R_+,H_0)\cap H_0(\R_+, H_k)\longrightarrow
   H_k(\R,H_0)\cap H_0(\R,H_k).
\]
The assertion follows from complex interpolation and \myref{I.G2-2.22}.
\end{proof}
\begin{remark} The construction of the extension operator $E^{(N)}$
used here is standard by now (cf. e.g. \cite[Sec. 4.4]{Tay:PDEI}). 
Seeley \cite{See:EFDHS} refined the construction of the sequence
$(a_j)$ in \myref{I.G2-2.26} and showed that $E^{(N)}$ can indeed
be constructed in such a way that it is independent of $N$. Such an extension
operator extends $C^\infty$--functions defined in a half
space to $C^\infty$--functions in the whole space. 
As remarked by Seeley the method of construction
of $E^{(N)}$ can be traced back to work of L. Lichtenstein
\cite{Lic:EBRA} and M.R. Hestenes \cite{Hes:ERDF}.
\end{remark}
\begin{cor}\label{I.S2-2.16} 
\begin{thmenum}
\item $({\ch}_s(\R_+, A))_{s\in\R_+}$ is a scale of Hilbert spaces.
\item For $s\in\R_+$ we have
\[
 {\ch}_s(\R_+, A) = H_s(\R_+,H_0)\cap H_0(\R_+,H_s)=
    \bigcap_{0\le t\le s} H_t(\R_+,H_{s-t}).
\]
The norm estimates \myref{I.G6-2.24} hold with $\R_+$ in place
of $\R$.

Thus, the operators $E^{(N)}$ defined in Proposition 
\plref{I.S2-2.15} have unique continuous extensions
\[
  E^{(N)}_s :{\ch}_s(\R_+,A) 
  \to {\ch}_s (\R, A), \quad s\le N-1,
\] 
satisfying $(E^{(N)}f) \upharpoonright \R_+=f.$
\end{thmenum}
\end{cor}
\begin{proof} (1) Observe first that in view of \myref{I.G2-2.12}
there is a continuous inclusion map
\begin{equation}
{\ch}_s (\R_+, A) \subset H_s(\R_+,H_0)\cap H_0(\R_+,H_s).
   \label{I.G2-2.27}
\end{equation}

Denote by $R(f):= f\upharpoonright \R_+ $ the restriction map onto
$\R_+$. By definition we have for $s<t$, $\theta\in [0,1]$
\begin{align*}
 {\ch}_{\theta t+(1-\theta)s}(\R_+, A) &=R ({\ch}_{\theta t+(1-\theta)s}(\R,A))\\
 &= R ([{\ch}_s(\R,A), {\ch}_t (\R,A)]_\theta )\\
 &\subset [{\ch}_s(\R_+,A), {\ch}_t(\R_+,A)]_\theta .
\end{align*}
To prove the converse let $N\ge t+1$ be an integer.
Then, in view of \myref{I.G2-2.27} 
and Proposition \plref{I.S2-2.15},
for 
$f\in[{\ch}_s(\R_+,A), {\ch}_t (\R_+,A)]_\theta 
%\cup (H_{\theta t+(1-\theta)s} (\R_+,H_0)
%\cap H_0(\R_+,H_{\theta t+(1-\theta)s}))
$
we have
$ E^{(N)}(f) \in
[{\ch}_s(\R,A), {\ch}_t(\R,A)]_\theta=
 {\ch}_{\theta t+(1-\theta)s}(\R,A)$ 
and thus $f=R(E^{(N)}f)\in {\ch}_{\theta t+(1-\theta)s}(\R_+,A).$
This proves (1).

(2) It suffices to prove the converse inclusion of \myref{I.G2-2.27}.
Let $f\in H_{s} (\R_+,H_0)\cap H_0(\R_+,H_{s})$. If $N\ge s+1$ is
an integer then, as before, we have  $E^{(N)}f\in 
\ch_s(\R,A)$ and thus
$f=R(E^{(N)}f)\in {\ch}_{s}(\R_+,A).$

The norm estimates follow from the fact that
two comparable Hilbert space norms are equivalent
and from the continuity of the inclusion \myref{I.G2-2.27}.
\end{proof}

We denote by $\cl^p(H)$ the von Neumann--Schatten class of $p$--summable
operators in $H$. A linear map $T:H\to \tilde H$ from the Hilbert space
$H$ into the Hilbert space $\tilde H$ is in the class $\cl^p$
if $T^*T\in \cl^{p/2}(H)$; this implies $TT^*\in\cl^{p/2}(\tilde H)$.

\begin{prop}\label{I.S1-2.13} If $(A+i)^{-1}$ is in $\cl^p(H)$  
(compact), 
then for any $\varphi\in\cinfz{\R}$ the map
$\ch_1(\R_+,A)\hookrightarrow
L^2(\R,H), f\mapsto \varphi f$, is of class $\cl^{p+1}$
(compact).
\end{prop}

\begin{proof} W.l.o.g. we may assume $\supp\varphi\subset (0,1)$. We consider
the operator $\tau_0=-\frac{d^2}{dx^2}+A^2$ on
%\begin{equation}
\[\{f\in\cinf{[0,1],H_\infty}\,|\, f(0)=f(1)=0\}.\]
%\end{equation}
Since $(A+i)^{-1}$ is compact, the operator $A$ is discrete. By slight
abuse of notation let $(e_a)_{a\in\spec A}$ be an orthonormal
basis of $H$ with $A e_a=a e_a$. Clearly, $e_a\in H_\infty$.
Then
\begin{equation}
   \{\sqrt{2}\sin(n\pi .)\otimes e_a\}_{n\in\N,a\in A}
  \label{I.G1-2.15}
\end{equation}
is an orthonormal basis of $L^2([0,1],H)$. Furthermore,
\begin{equation}
   \tau_0(\sqrt{2}\sin(n\pi .)\otimes e_a)=(n^2\pi^2+a^2)
    \sqrt{2}\sin(n\pi .)\otimes e_a.
   \label{I.G1-2.16}
\end{equation}
Thus \myref{I.G1-2.15} is an orthonormal basis of eigenvectors of $\tau_0$.
Hence $\tau_0$ is essentially self--adjoint and $\ovl{\tau}_0$ is
discrete. This implies that the map $\cd(\ovl{\tau}_0)\hookrightarrow
L^2([0,1],H)$ is compact.

If $(A+i)^{-1}\in \cl^p(H)$ then
\begin{align*}
    \sum_{n=1}^\infty &\sum_{a\in\spec A} (1+ n^2\pi^2 +a^2)^{-p/2-1/2}\\
 &\le C \sum_{a\in\spec A} \int_0^\infty (a^2+t^2+1)^{-p/2-1/2}dt\\
 &\le C' \sum_{a\in\spec A} (a^2+1)^{-p/2}<\infty,
\end{align*}
hence the map $\cd(\ovl{\tau}_0)\hookrightarrow L^2([0,1],H)$
is of class $\cl^{(p+1)/2}$. Since the map
\[\ch_2(\R_+,A)\hookrightarrow L^2(\R,H), f\longmapsto \varphi f,\]
factorizes through
$\cd(\ovl{\tau}_0)$ we conclude that 
$\ch_2(\R_+,A)\hookrightarrow L^2(\R,H), 
f\mapsto \varphi f$, is of class  $\cl^{(p+1)/2}$, too.
Now the assertion follows from interpolation.\end{proof}

\section{Fredholm pairs}
\label{Secthree}

\subsection{Fredholm pairs in a Hilbert space}
\newcommand{\PH}{\cp(H)}
Let $H$ be a Hilbert space. We denote by $\PH$ the set of orthogonal
projections on $H$.

\begin{dfn}\label{I.S1-3.1} Let
$P,Q\in\PH$. The pair $(P,Q)$  
is called a \textnormal{Fredholm pair} if $Q:\im P\to \im Q$ is
a Fredholm operator. The index of this operator is denoted by
$\ind(P,Q)$.

The pair $(P,Q)$ will be called \textnormal{invertible} if $Q:\im P\to \im Q$
is invertible. 
\end{dfn}

We will see below that these definitions are symmetric in $P$ and $Q$.
The notion of a Fredholm pair 
was introduced by Kato \cite[IV.4.1]{Kat:PTLO}. Bojarski
\cite{Boj:ALCPFPS} seems to be the first one who used this concept
in the theory of elliptic boundary value problems (cf. also
Booss--Wojciechowski \cite[Sec. 24]{BooWoj:EBPDO},\cite{BooMorStrWoj:GCA}).
Recently, the notion was systematically studied by
Avron, Seiler, and Simon \cite{\ASS} who apparently were not aware
of the earlier literature.
The following fact is proved by straightforward calculation
\cite[Sec. 2]{\ASS}:

\begin{prop}\label{I.S1-3.2}
Let $P,Q\in \PH$. We put
  \begin{equation}
  \begin{split}
  X&:= P-Q,\quad Y:=I-P-Q,\\
  K&:= 2PQ-P-Q=-(2P-I)X=X(2Q-I).
  \end{split}
   \label{I.G1-3.1}
  \end{equation}  
Then 
\begin{thmenum}
\item $X^2+Y^2=I,$ $XY=-YX$, and $X^2$ commutes with $P$ and $Q$.
\item $KK^*=K^*K=X^2$, in particular, $K$ is a normal operator. Moreover,
\[\|K\|\le \|P-Q\|.\]
\item $Z:=I+K$ satisfies $ZQ=PZ=PQ$.
\end{thmenum}
\end{prop}

The following proposition gives a
useful Fredholm criterion.

\begin{prop}[{\cite[Prop. 3.1]{\ASS}}]\label{I.S1-3.3}
$(P,Q)$ is a Fredholm pair if and only if
\begin{equation}
    \pm 1\not\in \specess(P-Q).
    \label{I.G1-3.2}
\end{equation}
In this case,
\begin{equation}
   \begin{split}
    \ker Q\cap\im P&=\ker (P-Q-I),\\
    \im Q\cap\ker P&=\ker (P-Q+I),
    \end{split}
    \label{I.G1-3.3}
\end{equation}
in particular,
\begin{equation}
   \ind(P,Q)=\dim\ker(P-Q-I)-\dim\ker(P-Q+I).
   \label{I.G1-3.4}
\end{equation}
\end{prop}
We see from \myref{I.G1-3.2} that, indeed,
$(P,Q)$ is a Fredholm pair if and only
if $(Q,P)$ is.

\begin{cor}\label{I.S1-3.4} For $P,Q\in\PH$ the following statements are
equivalent:
\begin{thmenum}
\item $\|P-Q\|<1$,
\item $\pm 1\not\in\spec(P-Q)$,
\item $(P,Q)$ is an invertible pair,
\item $l(PQ)=P,\, r(PQ)=Q$ and $PQ$ has closed range (where $l$ and $r$
denote the left and right support, respectively).
\end{thmenum}
\end{cor}
\begin{proof}  The equivalence
of (1) and (2) follows since 
$P-Q$ is self--adjoint and $\|P-Q\|\le 1$.
The equivalence of (2) and (3) follows from 
Proposition \plref{I.S1-3.3} and the equivalence of (3) and (4)
is well--known.\end{proof}

Since (1) is symmetric in $P,Q$ this corollary 
implies that if $(P,Q)$
is invertible then $(Q,P)$ is invertible, too.

\begin{remark} \label{I.S1-3.9} 
We note for future reference that $(P,Q)$ is Fredholm if
and only if the operator
\begin{equation}
T:=PQ+(I-P)(I-Q)\label{I.G5-3.5}
\end{equation}
is a Fredholm operator. In this case
\begin{equation}
  \ker T = \ker T^*= \im Q \cap \ker P \oplus \ker Q \cap \im P,
\end{equation}
hence $(P,Q)$ is an invertible pair if and only if $T$ is invertible.
Note that $\ind T$ is always zero.
\end{remark}

We will now study deformations of Fredholm pairs.

\begin{lemma}\label{I.S1-3.5} Let $P,Q\in\PH$ with $\|P-Q\|<1$.
Then $Z_s:=I+sK, 0\le s\le 1$, with $K$ from \myref{I.G1-3.1},
 is an invertible and normal
operator, $P_s:=Z_sQZ_s^{-1}\in \PH$ and
\[P_1=P,\quad \|P_s-Q\|<1,\quad 0\le s\le 1.\]
\end{lemma}
\begin{proof} From Proposition \plref{I.S1-3.2} we infer that $K$ is normal
with $\|K\|<1$,
hence $Z_s$ is normal and invertible. Furthermore, 
$P_1=ZQZ^{-1}=P$. 
Moreover,
\begin{equation}
  Z_s^*Z_s=I+s(K+K^*)+s^2X^2=I+(s^2-s)X^2
  \label{I.G1-3.7}
\end{equation}
commutes with $P$ and $Q$.
Since $Z_s$ is normal we have $P_s\in\PH$. We thus also have
\begin{equation}
    P_s=U_sQU_s^*
    \label{I.G1-3.8}
\end{equation}
with the unitary operator $U_s=Z_s(Z_sZ_s^*)^{-1/2}$.

Since $\|Z_s-I\|\le \|K\|<1$ and since $Z_s$ is normal, the
spectrum of $U_s$ is contained in $\{z\in\C\,|\,|z|=1, 
|\operatorname{arg}\,z|\le \pi/2-\delta\}$ for some $\delta>0$. 
Hence there exists
$\eps>0, 0<q<1$ such that for $z\in \spec U_s$ we have
\[|z-\eps|\le q.\]
Consequently,
\begin{equation}
   \|P_s-Q\|=\|[U_s-\eps I,Q]\|\le q<1.
   \label{I.G1-3.9}
\end{equation}
The last estimate uses the fact that
for $P\in \PH$ and any bounded operator $T\in \cl(H)$
we have the estimate
\begin{equation}
  \|[P,T]\|\le \|T\|.
   \label{I.G1-3.5}
\end{equation}
This follows immediately from
\begin{equation}
   [P,T]=PT-TP=PT(I-P)-(I-P)TP.
   \label{I.G1-3.6}
\end{equation}
The lemma is proved.\end{proof}

Next we consider an arbitrary Fredholm pair $(P,Q)$, $P,Q\in\PH$.
We put
\begin{equation}\begin{aligned}[c]
      P'&:=l(PQ),\\
      P''&:=P-P',
  		\end{aligned}\quad
     \begin{aligned}[c]
        Q'&:=r(PQ),\\
        Q''&=Q-Q''.
     \end{aligned}
   \label{I.G1-3.10}
\end{equation}
Note that $P''$ is the orthogonal projection onto $\im P\cap \ker Q$
and $Q''$ is the orthogonal projection onto $\ker P\cap \im Q$.

Since $PQ=P'Q'$ and since $PQ$ has closed range we infer from
Corollary \plref{I.S1-3.4} that $(P',Q')$  is an invertible pair,
hence

\begin{equation}
    \|P'-Q'\|<1.\label{I.G1-3.11}
\end{equation}

Furthermore, by construction
\begin{equation}
     P'\le P,\quad Q'\le Q.
   \label{I.G1-3.12}
\end{equation}

Let $U_s=U_s(P',Q')$ be the family of unitaries constructed in the proof
of Lemma \plref{I.S1-3.5}. Then
\begin{equation}
      P'=U_1Q'U_1^*.
   \label{I.G1-3.13}
\end{equation}
Put
\begin{align}
    \tilde Q&:=U_1^*PU_1=Q'+U_1^*P''U_1,\label{I.G1-3.14}\\
\intertext{and}
     P_s&:=U_s\tilde Q U_s^*,\quad  P_s':=U_s Q'U_s^*,\quad 0\le s\le 1.
    \label{I.G1-3.15}
\end{align}
Then 
one sees that $\|P_s-\tilde Q\|<1$, as in the proof Lemma \plref{I.S1-3.5},
and
since $\tilde Q-Q'$ is of finite rank, we see that
$(P_s,Q)$, $0\le s\le 1$, is
a continuous family of Fredholm pairs with
\begin{equation}
    P_0=\tilde Q=Q'+\tilde Q'',\quad P_1=P.
    \label{I.G1-3.16}
\end{equation}

Thus we have proved the following fact.

\begin{lemma}\label{I.S1-3.6} Let $(P,Q)\in \PH$ be a Fredholm pair.
Then there is a continuous family of Fredholm pairs $(P_s,Q)$,
$0\le s\le 1$, such that $P_1=P$ and $P_0=\tilde Q=Q'+\tilde Q''$
with a finite rank projection $\tilde Q''\le (I-Q')$.
More precisely, $P_s=U_s\tilde QU_s^*$ where $U_s$ is
constructed from $P',Q'$ as in the proof of Lemma \plref{I.S1-3.5}.
\end{lemma}

Consider the Fredholm pair $(\tilde Q,Q)$ constructed above.
We  have
\begin{equation}
   \begin{split}
      \dim(\ker\tilde Q\cap\im Q)&=\dim (\ker \tilde Q''\cap \im Q'')\\
                                 &=\rank( Q'')-\rank (\tilde Q''Q'')\\
      \dim(\im\tilde Q\cap\ker Q)&=\dim (\im \tilde Q''\cap \ker Q'')\\
                                 &=\rank (\tilde Q'')-\rank (\tilde Q''Q''),
        \end{split}
    \label{I.G1-3.17}
\end{equation}
hence
\begin{equation}
   \ind(P,Q)=\rank\tilde Q''-\rank Q''.
 \label{I.G1-3.18}
\end{equation}
If $\ind(P,Q)=0$ we can find a path of orthogonal projections of finite rank
in the space $(I-Q')(H)$
connecting
$\tilde Q''$ and $Q''$, and hence we can find a path $(P_s,Q)$ of 
Fredholm pairs connecting $(P,Q)$ with
the pair $(Q,Q)$.

\begin{theorem}\label{I.S1-3.7} For fixed $Q\in\PH$ the connected components
of 
\[\bigsetdef{P\in\PH}{(P,Q)\;\text{Fredholm}\;}\]
are labeled by $\ind(P,Q)$.

More precisely, given $P,P'\in \PH$ such that the pairs $(P,Q),(P',Q)$
are Fredholm with the same index
then there is a smooth family of unitary operators $U_t,
0\le t\le 1$, such that
\begin{nonumabsatz}
\item $U_0=I$, $U_1PU_1^{-1}=P'$
\item $(U_tPU_t^{-1},Q)$ is Fredholm.
\end{nonumabsatz}
\end{theorem}

\mathversion{bold}
\subsection{Fredholm pairs of $\gamma$--symmetric projections of order $0$}
\mathversion{normal}

In the discussion of self-adjoint extensions of 
$D=\gamma(\frac{d}{d x} + A)$ below we will need 
projections $P\in {\cp} (H)$ which are 
$\gamma$-\emph{symmetric} in the sense that
\begin{equation}
  \gamma P = (I - P)\gamma.
  \label{I.G2-3.19}
\end{equation}
In particular, we will require the existence of a 
specific projection with \myref{I.G2-3.19}: 

\medskip\noindent
\begin{assumption} \label{I.S2-3.8}
There is a spectral projection, $P_+(A),$ of $A$ satisfying
\myref{I.G2-3.19} and, in addition,
\begin{equation}
  1_{[0,\infty)}(A) \ge P_+(A) \ge 1_{(0,\infty)}(A).
   \label{I.G2-3.20}
\end{equation}
\end{assumption}
We will write $P_+:=P_+(A)$ if no confusion is possible.
Furthermore we abbreviate
\begin{equation}
\Pmi:=\Pmi(A):=I-\Ppl(A),\label{I.G3-4.19}
\end{equation}
In view of our basic structural assumptions we easily see the following
fact.

\begin{prop}\label{I.S2-3.9}
A projection with \myref{I.G2-3.19} and \myref{I.G2-3.20} 
exists if and only if  the involution  
$i\gamma \,|\, \ker A$ has signature $0$ 
by which we mean that $\ker(\gamma+ i)\cap \ker A\simeq \ker(\gamma-i)\cap A$.
A unique such projection exists if and only if $\ker A=0$.
\end{prop}

Namely, if $i\gamma \,|\, \ker A$ has signature $0$
we may choose an isometry 
\begin{equation}
U:\ker(\gamma+ i)\cap \ker A
\longrightarrow \ker(\gamma-i)\cap A
 \label{I.G7-3.24}
\end{equation}
and put
\begin{equation}\sigma:\ker A\longrightarrow \ker A,\quad 
\sigma:=\begin{pmatrix} 0 & U\\ U^*& 0 \end{pmatrix}.
\end{equation}
Then the orthogonal projection
\begin{equation}
P_+(A,\sigma):= 1_{(0,\infty)}(A)\oplus \frac{I+\sigma}{2}
\end{equation}
satisfies Assumption \plref{I.S2-3.8} and it is clear that all
such projections are in one--one correspondence with unitaries
$U:\ker(\gamma+ i)\cap \ker A\longrightarrow \ker(\gamma-i)\cap A$.

If $i\gamma$ has signature $\not= 0$ on $\ker A$, we can often remedy this 
defect by slightly extending $H$. Namely, w.l.o.g. let $V$ be a Hilbert space
with
\begin{equation}
      \ker(\gamma+ i)\cap \ker A\simeq 
        \bigl(\ker(\gamma-i)\cap A\bigr)\oplus V.
      \label{I.G5-3.24}
\end{equation}
If $\dim\ker A<\infty$ then $V$ is finite--dimensional, too. We 
put
\begin{equation}
  \tilde H:= H\oplus V,\quad
  \tilde A:=A\oplus 0,\quad
  \tilde \gamma:= \gamma\oplus i.
\end{equation}
In view of \myref{I.G5-3.24} we have 
$\ker (\tilde\gamma+i)\cap\ker \tilde A)\simeq
\ker (\tilde \gamma-i)\cap\ker \tilde A$. Hence there exists 
a unitary operator $\sigma:
\ker\tilde A\to \ker \tilde A$ with $\sigma^2=I, \sigma\tilde \gamma=-
\sigma\tilde\gamma$, and thus the orthogonal projection
\begin{equation}
   P_+(\tilde A):=1_{(0,\infty)}(A)\oplus\frac{I+\sigma}{2}
\end{equation}
satisfies Assumption \plref{I.S2-3.8} with respect to 
$\tilde A, \tilde \gamma$.

\begin{dfn}\label{I.S2-3.11}
Let ${\ca}\subset \OpAn$ be an operator algebra. 
We introduce
\[ %\begin{split}
{\cp}_s({\ca}) :=\biggl \{ 
P\in {\cp} (H)\cap {\ca} \,\biggm|\, \begin{array}{l}
    (P, P_+) \,\text{is a Fredholm pair and}\\ 
    P\, \text{is $\gamma$-symmetric}
   				     \end{array}\biggr\}.
%   \end{split}
\]
\end{dfn}
Note that for $P\in\cp_s(\ca)$ we have
$\gamma(P-\Ppl)=(\Ppl-P)\gamma$, hence
\begin{equation}
   \ind(P,\Ppl)=0.
     \label{I.G1-3.19}
\end{equation}

We want to derive the analogue of Theorem \plref{I.S1-3.7} for 
${\cp}_s ({\ca})$. 
For this to work we need an additional structural property of 
the algebra in question. 
Thus from now on we consider algebras, $\PsiA\subset \OpAn$, 
with the following properties (cf. \myref{axiomsintro} in the introduction): 

\label{axiomspaperbody}
%\begin{itemize}
%\item[($\Psi$1)] $\Psi^0(A)$ is a $\ast$-subalgebra of
%$\OpAn$ with \emph{holomorphic functional calculus};\\
%\item[($\Psi$2)] $\Psi^0(A)$ contains an orthogonal projection,
%$P_+(A)$, satisfying 
%\begin{equation}
% I-P_+(A)=\gamma^*P_+(A)
%\gamma,\quad P_{(0,\infty)} (A)\le P_+(A)\le P_{[0,\infty)}(A).
% \label{intro-1.19}
%\end{equation}
%\end{itemize}

%\begin{alphnumbering}[I.G3-3.22]
%\begin{MLlist}
%\item $\PsiA$ contains $P_+$ and all finite rank operators in $\OpAn$,
%      \puteqnnum\label{I.G3-3.22a}
%\item $\PsiA$ admits holomorphic functional calculus.
%      \puteqnnum\label{I.G3-3.22b}
%\end{MLlist}
%\end{alphnumbering}
\begin{alphnumbering}[I.G3-3.22]
\begin{MLlist}
\begin{MLnumlist}
\MLitem{$\PsiA$ is a $\ast$-subalgebra of
$\OpAn$ containing the smoothing operators and
with holomorphic functional calculus;}
      \puteqnnum\label{I.G3-3.22b}

\medskip
\MLitem{$\PsiA$ contains an orthogonal projection $P_+$ with the properties
stated in Assumption \plref{I.S2-3.8}.} \puteqnnum\label{I.G3-3.22a}
\end{MLnumlist}
\end{MLlist}
\end{alphnumbering}

Assuming \myref{I.G3-3.22b},
in view of Proposition \plref{I.S2-3.9}
property \myref{I.G3-3.22a} is satisfied if we can find
$U \in \Psi^0 (A)$ as in \myref{I.G7-3.24}, 
since all finite spectral projections of $A$ are smoothing.

Then the above construction goes through to yield the following result.

\begin{theorem}\label{I.S1-3.10} $\cp_s(\PsiA)$ is path connected. More precisely,
given $P,Q\in\cp_s(\PsiA)$ there is a path 
$P_t\in \cp_s(\PsiA), 0\le t\le 1$,
connecting $P$ and $Q$ such that $t\mapsto P_t$ is smooth for all
$H_s$--norms.
$P_t$ can be chosen of the form $P_t=U_tP_0U_t^{-1}$
where $U_t$, for $0\le t\le 1$, is a smooth family of unitaries
in $\PsiA$ satisfying
\begin{MLlist}
\item $U_0=I, U_t\gamma=\gamma U_t,$
\item and $t\mapsto U_t$ is smooth for all $H_s$--norms.
\end{MLlist}
\end{theorem}
\begin{proof}
Let $X=P-\Ppl, Y=I-P-\Ppl, K=2P\Ppl-P-\Ppl$ (cf. Proposition \plref{I.S1-3.2}).
Note the relations
\begin{equation}\begin{aligned}[c]
      \gamma X&=-X\gamma,\\
       \gamma K&=K\gamma,
     		\end{aligned}\quad
        \begin{aligned}[c]
           &\gamma Y=-Y\gamma\\
           &X,Y,K\in\PsiA.
 	\end{aligned}		
   \label{I.G1-3.20}
\end{equation}
Furthermore, $Z_s=I+sK\in\PsiA$ and commutes with $\gamma$.
Thus if $\|P-\Ppl\|<1$ then in view of Lemma \plref{I.S1-3.5} and 
\myref{I.G1-3.8}
\begin{equation}
   P_s=Z_s\Ppl Z_s^{-1}=U_s\Ppl U_s^*,\quad U_s=Z_s(Z_sZ_s^*)^{-1/2}
     \label{I.G1-3.21}
\end{equation}
is the desired family of projections in $\cp_s(\PsiA)$ connecting $P$ and
$\Ppl$. Here, we only used the spectral
invariance of $\PsiA$.

Next we consider a general $P\in \cp_s(\PsiA)$. To show that there is
a path from $P$ to $\Ppl$ having the desired properties we 
repeat the construction
\myref{I.G1-3.10}--\myref{I.G1-3.16}. The invariance under holomorphic functional calculus of
$\PsiA$ implies $P''\in\PsiA$ and thus $P'\in\PsiA$.
To see this note that $P''$ is the orthogonal projection
onto $\im P\cap \ker \Ppl=\ker(P-\Ppl-I)$ and $0\not\in 
\specess(P-\Ppl-I)$  by Proposition \plref{I.S1-3.3}.
Thus $P''$ is an analytic function of $P-\Ppl-I\in\PsiA$.
Analogously, we have $\Ppl', \Ppl''\in\PsiA$.

However, since $P'\le P$, it is $\gamma$--symmetric only if $P''=0$,
hence $P'\not\in\cp_s(\PsiA)$ if $P''\not=0$.
We next show that $P'+P_+''\in\cp_s(\PsiA).$

$\Ppl''$ is the orthogonal projection onto $\im \Ppl\cap \ker P=
\ker (I-\Ppl)\cap\im (I-P)$, thus by symmetry
\begin{equation}
    \gamma P''\gamma^*=\Ppl''.
    \label{I.G1-3.22} 
\end{equation}

Since $\im \Ppl''\perp\im P'$, $P'+\Ppl''\in\PsiA$ is an orthogonal projection
and, by \myref{I.G1-3.22}, it is $\gamma$--symmetric.
Furthermore,
by \myref{I.G1-3.11}
\begin{equation}
    \|(P'+\Ppl'')-\Ppl\|=\|P'-\Ppl'\|<1,
   \label{I.G1-3.23}
\end{equation}
thus $(P'+\Ppl'',\Ppl)$ is an invertible pair.
Thus by the first part of this proof, there is a path from
$P'+\Ppl''$ to $\Ppl$ in $\cp_s(\PsiA)$ with the desired properties.

As similar argument as in the proof of Theorem \plref{I.S1-3.7}
now shows that there also is such a path from $P'+\Ppl''$ to $P$.
\end{proof}

\section{Regularity for the model operator}
\label{Secfour}

In this section we study the operator $D=D_+$ and the associated boundary
value problems on $\R_+$.

\newcommand{\dpl}{D}

In view of \myref{I.G1-4.1} it is natural to consider, for an
orthogonal projection
$P\in\OpAn$, the operator $D_P$  given by 
\begin{equation}\begin{split}
     \cd(D_P)&:=\bigsetdef{u\in \ch_1(\R_+,A)}{Pu(0)=0},\\
     D_P&:=D^*\restr\cd(D_P).
   		\end{split}
   \label{I.G1-4.2}
\end{equation}
Since $P$ is of order zero and, by Lemma \plref{I.S1-4.1},
$u(0)\in H_{-1/2}$ for $u\in\cd(D^*)$ we have a natural extension
of $D_P$ (which we denote by $D_{P,\max}$) given by
\begin{equation}\begin{split}
   \cd(D_{P,\max})&:=\bigsetdef{u\in\cd (D^*)}{Pu(0)=0},\\
   D_{P,\max} &:= D^*\restr \cd(D_{P,\max}).
  		\end{split}
   \label{I.G1-4.4}
\end{equation}
We are mainly interested in those projections which render $D_P$
self--adjoint. To characterize them, we prepare two results.

\begin{prop}\label{I.S1-4.2}
\begin{thmenum}
\item $D_P$ is symmetric if and only if 
$I-P\le \gamma P\gamma^*$.
\item Let  $P_1, P_2$ be orthogonal projections in $\OpAn$.
Then $P_1\le P_2$ $(P_1\not=P_2)$ implies $D_{P_2}\subset D_{P_1}$
$(D_{P_1}\not= D_{P_2})$.
%\item If $(I-P)<\gamma P\gamma^*$ then there exists a projection
%$P_1$ of order $0$ with $(I-P_1)\le\gamma P_1\gamma^*$ and
%$P_1<P$.
\end{thmenum}
%Thus if $(I-P)< \gamma P\gamma^*$ then $D_P$ has a proper symmetric
%extension.
\end{prop}
\begin{proof} 
(1) By Lemma \plref{I.S1-4.1}, we have for $u,v\in \cd(D_P)$ 
\[(-D_Pu,v)+(u,D_Pv)=\scalar{\gamma u(0)}{v(0)}.\]
Thus $D_P$ is symmetric if and only if 
for all $\xi, \eta\in\ker\,(P\restr H_{1/2})$
\[\scalar{\gamma \xi}{\eta}=0,\]
the "only if " part following from Proposition \plref{I.S1-2.11}.
But this is easily seen to be equivalent to 
$I-P\le \gamma P\gamma^*$.

(2) $P_1\le P_2$ implies
$D_{P_2}\subset D_{P_1}$. If $P_1\le P_2, P_1\not=P_2,$ then we choose
$\xi\in (H_{1/2}\setminus\{0\})\cap (\ker P_2\setminus \ker P_1)$. 
By Proposition \plref{I.S1-2.11}, we can find $f\in \ch_1(\R_+,A)$
with $f(0)=\xi$. Then $f\in \cd(D_{P_2})\setminus \cd(D_{P_1})$.
\end{proof}

%3. Choose $\xi\in (H_\infty\setminus\{0\})\cap\im(\gamma P\gamma^*-(I-P))$
%with $\gamma \xi\perp\xi$. This is possible since
%$\im(\gamma P\gamma^*-(I-P))$ is $\gamma$--invariant.
%
%Let $I-P_1$ be the orthogonal projection onto $\im(I-P)\oplus <\xi>$.
%Then $P_1$ does the job.\end{proof}

We will show in Theorem \plref{I.S1-4.4} below that 
if $(I-P)\lvertneqq\gamma^*P\gamma$ then $D_P$ is not
self--adjoint.

We abbreviate for $a\in\R$
\begin{equation}\begin{aligned}[c]
   P_{<a}&:=1_{(-\infty,a)}(A),\\
    P_{>a}&:=I-P_{\le a},
   		\end{aligned}\quad
     \begin{aligned}[c]
   P_{\le a}&:=1_{(-\infty,a]}(A),\\
   P_{\ge a}&:=I-P_{<a}.
     \end{aligned}
\end{equation}

\begin{lemma}\label{I.S1-4.3} Let $r^*:\cd(D^*)\to H_{-1/2}$ be the
restriction map from Lemma \plref{I.S1-4.1}.
We have
\begin{alphnumbering}[I.G1-4.5] 
\begin{equation}
    r^*(\cd(D^*))=P_{>0} (H_{-1/2})\oplus P_{\le 0} (H_{1/2}),
    \label{I.G1-4.5a}
\end{equation}
in particular
\begin{equation}
r^*(\cd(D_{P_{>0},\max}))=P_{\le 0} (H_{1/2})\subset H_{1/2}.
  \label{I.G1-4.5b}
\end{equation}
\end{alphnumbering}
\end{lemma}
\begin{proof} 
We choose
$\chi\in\cinfz{\R}, \chi = 1$ near $0$, and put
for $f\in\cd(D_{P_{>0,\max}})$
\begin{equation}
g(x):=   \chi(x) e^{-Ax}f(0), \quad x\le 0.
   \label{I.G1-4.6}
\end{equation}
Then $g\in\cd(D_-)^*, g(0)=f(0)$, so from Proposition \plref{I.S2-2.13}
we infer $f\in\ch_1(\R_+,A)$ and hence $f(0)\in H_{1/2}$.
This proves $r^*(\cd(D_{P_{>0,\max}}))\subset H_{1/2}$.

To prove \myref{I.G1-4.5a} let $f\in \cd(D^*)$. 
Then by Lemma \plref{I.S1-4.1} we have
$f(0)\in H_{-1/2}$. Moreover, $x\mapsto P_{\le 0} f(x)$ lies in
$\cd(D_{P_{>0},\max})$ and hence by the first part of the proof
we have $P_{\le 0} f(0)\in H_{1/2}$.

Conversely, let $\xi\in P_{>0} (H_{-1/2}), \eta\in P_{\le 0} (H_{1/2})$ 
be given.
By Proposition \plref{I.S1-2.11} there exists $f\in\ch_1(\R_+,A)$
with $f(0)=\eta$. Since $P_{\le 0}\xi=0$ we see that \myref{I.G1-4.6},
with $\xi$ in place of $f(0)$, defines, for $x>0$, a function
$g$ in $\cd(D^*)$ with  $g(0)=\xi$.
Thus $\xi+\eta=r^*(f+g)$ and we reach the conclusion.
\end{proof}

Now we present the main result of this section.

\begin{theorem}\label{I.S1-4.4} 
Let $P\in\OpAn$ be an orthogonal projection.

\begin{alphnumbering}[I.G3-4.5]
\begin{thmenum}
\item $D_P^*=D_{\gamma^* (I-P)\gamma,\max}$.
\item The following statements are equivalent:
      \begin{thmenum}
         \item $D_P=D_{P,\max}$.\puteqnnum
         \item $r^*(\cd(D_{P,\max}))\subset H_{1/2}$.\puteqnnum\label{I.G3-4.5b}
         \item If $\xi\in H_{-1/2}$ satisfies $P\xi=0$ and 
   $P_{<0} \xi\in H_{1/2}$,\\ then also $\xi\in H_{1/2}$.\puteqnnum
   \label{I.G3-4.5c}
      \end{thmenum}
\end{thmenum}
\end{alphnumbering}
\end{theorem}
%From now on we {\em assume} that there exists a $\gamma$--symmetric
%orthogonal projection $\Ppl\in\OpAn$ with
%$1_{(0,\infty)}(A)\le \Ppl\le 1_{[0,\infty)(A)}$. We put $\Pmi:=I-\Ppl$.
%(cf. ($\Psi$2)).
%\item $\tilde D_P$ is self--adjoint.
%\end{enumerate}
%\end{prop}

\begin{proof} (1) For $f\in\cd(D_{\gamma^*(I-P)\gamma,\max})$ and $g\in \cd(D_P)$
we have by Lemma \plref{I.S1-4.1} (2)
\begin{equation}\begin{split}
    -(D^*f,g)+(f,D^*g)&=\Bss{1/2}{\gamma f(0)}{g(0)}\\
   &=\Bss{1/2}{P\gamma f(0)}{(I-P)g(0)}=0,
   		\end{split}
    \label{I.G3-4.7}
\end{equation}
since the last written sesquilinear form vanishes identically on
$H_\infty\times H_\infty$ and extends to $H_{-1/2}\times H_{1/2}$ by
continuity.
Hence $D_{\gamma^*(I-P)\gamma,\max}\subset D_P^*.$

To prove $D_P^*\subset D_{\gamma^*(I-P)\gamma,\max}$ 
we consider $f\in\cd(D_P^*)$.
Since $D_P^*\subset D^*$, \myref{I.G3-4.7} 
gives for all $g\in\cd(D_P)$
\begin{equation}
     0=\Bss{1/2}{\gamma f(0)}{g(0)}=\Bss{1/2}{(I-P)\gamma f(0)}{g(0)}.
     \label{I.S3-4.8}
\end{equation}
By Corollary \plref{I.S1-2.12}, this implies 
$\gamma^*(I-P)\gamma f(0)=0$ and hence 
\[f\in\cd(D_{\gamma^*(I-P)\gamma,\max}).\]

(2) (i)$\Rightarrow$(ii): This follows immediately from
the definition of $D_P$ and  
%If $D_P=D_{P,\max}$ then
%by 1. $\cd(D_{P,\max})=\cd(D_P)\subset \ch_1(\R_+,A)$ and 
%the assertion follows from 
the Trace Theorem
\plref{I.S1-2.10}.

(ii)$\Rightarrow$(i): 
In view of Corollary \plref{I.S2-2.14} (1), (ii) implies
that $\cd(D_{P,\max})\subset\ch_1(\R_+,A)$ and thus (i).

(ii)$\Rightarrow$(iii): 
Let $\xi\in H_{-1/2}, P\xi=0, P_{<0} \xi\in H_{1/2}$.
Then, in view of Lemma \plref{I.S1-4.3}, $\xi\in r^*(\cd(D^*))$; 
hence there exists
$f\in\cd(D^*)$ with $f(0)=\xi$. Since $P\xi=0$ 
we have $f\in\cd(D_{P,\max})$ and
thus $\xi=f(0)\in H_{1/2}$.

(iii)$\Rightarrow$(ii): If $f\in\cd(D_{P,\max})$ then, again
by Lemma \plref{I.S1-4.3}, we have $f(0)\in H_{-1/2},$ and
$P_{<0} f(0)\in H_{1/2}$; since
$Pf(0)=0$, (iii) implies $f(0)\in H_{1/2}$.
\end{proof}

Theorem \plref{I.S1-4.4} motivates the following terminology.

\begin{dfn}\label{I.S1-4.6} $P$ is called {\em regular} (with respect to $D$), 
if one of the equivalent conditions \myref{I.G3-4.5} is fulfilled.

Inductively, we call $P$ {\em $n$--regular} (with respect to $D$)
for  $n\ge 2$, if $P$ is $(n-1)$--regular
and
\[\cd(D_P^n)\subset \ch_n(\R_+,A).\]
\end{dfn}

Note that, for a regular projection, the restriction map 
$r^*:\cd(D_P)\to H_{1/2}$
is continuous. 

\begin{remark} One should note that $n$--regularity is a relative
notion, i.e. in view of \myref{I.G3-4.5c} it depends on $\Ppl$. More
precisely, we should therefore refer to the $n$--regular \emph{pair}
$(P,\Ppl)$.
However, since $A$ is fixed once and for all we avoid this
notion for simplicity. Referring to the pair $(P,\Ppl)$ 
cannot be avoided, however, in the case of ``Fredholmness''
(Definition \plref{I.S1-3.1})   and
``ellipticity'' (Definition \plref{I.S1-5.4}). 
The reason is that the pair
$(P,\Ppl)$ can be Fredholm (resp. elliptic) without $P$ having this
property.
\end{remark}

We note two consequences of Theorem \plref{I.S1-4.4}.

\begin{cor}\label{I.S3-4.4}Let $P\in\OpAn$ be an orthogonal projection.
Then $D_P$ is self--adjoint if and only if $P$ is regular and
$\gamma$--symmetric (in the sense of
\myref{I.G2-3.19}).
\end{cor}
\begin{cor}\label{I.S1-4.5} Let $P\in\OpAn$ be an
orthogonal projection.
Then
\[\cd:=\{u\in\cinfz{[0,\infty),H_\infty}\,|\, Pu(0)=0\}\]
is a core for $D_P$.
\end{cor}
\begin{proof} For the moment we put $T:=D_P\restr \cd$. Certainly, we have
$D\subset T\subset D_P$ and thus $D_P^*\subset T^*\subset D^*$.
We will show $T^*\subset D_P^*=D_{\gamma^*(I-P)\gamma,\max}$;
then $T^*=D_P^*$ and
hence $\ovl{T}=T^{**}=D_P^{**}=\ovl{D_P}$.

For $f\in \cd(T^*)$ and $g\in \cd(T)$ we find as in the previous
proof
\[ \Bss{1/2}{\gamma f(0)}{g(0)}=0.\]

Since $\ker \gamma^*(I-P)\gamma\cap H_\infty$ 
is dense in $\ker \gamma^*(I-P)\gamma$ this implies
$\gamma^*(I-P)\gamma f(0)=0$ and thus 
$f\in\cd(D_{\gamma^*(I-P)\gamma,\max})=\cd(D_P^*)$.
\end{proof}

We add a few comments on why the equivalent conditions
of Theorem \plref{I.S1-4.4} should be referred to as "regularity".

Consider the equation
\begin{equation}
      D^*f=g
     \label{I.G2-4.11}
\end{equation}
with $f,g\in L^2(\R_+,H)$. In general, \myref{I.G2-4.11} does not
imply $f\in \ch_1(\R_+,A)$, but Theorem \plref{I.S1-4.4} tells us that
\begin{equation}\begin{split}
      &D^*f=g,\quad f,g\in L^2(\R_+,H),\\
      & Pf(0)=0
		\end{split}
     \label{I.G2-4.12}
\end{equation}
implies $f\in\ch_1(\R_+,A)$ if and only if $P$ is regular.
$n$--regularity can be characterized similarly:

\begin{prop}\label{I.S2-4.7} An orthogonal projection $P\in\OpAn$ is $n$--regular
if and only if the relations
\begin{equation}
      \begin{split}
       &D^*f=g\in \ch_k(\R_+,A),\quad f\in L^2(\R_+,H),\\
       &Pf(0)=0,
      \end{split}
    \label{I.G2-4.13}
\end{equation}
imply $f\in \ch_{k+1}(\R_+,A)$, for $k\in\Z_+, 0\le k \le n-1$.
\end{prop}
\begin{proof} 
Let $P$ be a $n$--regular projection and consider the relations
\myref{I.G2-4.13} for some $k$.
Put 
\[ f_1:= f- e^{(k+1)}(0,(Df)(0),\ldots,(D^{k}f)(0)),\]
where $e^{(k+1)}$ is defined in Corollary \plref{I.S1-2.12}.

$f_1$ is well--defined since $Df\in\ch_{k}(\R_+,A)$. By construction,
$f_1\in\cd(D_P^{k+1})$, which is contained in 
$\ch_{k+1}(\R_+,A)$ since $P$ is $n$--regular.
But since $(D^jf)(0)\in H_{k+1/2-j}, 1\le j\le k$, we have
$ e^{(k+1)}(0,(Df)(0),\ldots,(D^{k}f)(0))\in \ch_{k+1}(\R_+,A)$ hence
$f\in\ch_{k+1}(\R_+,A)$.

Conversely, if \myref{I.G2-4.13} implies  $f\in \ch_{k+1}(\R_+,A)$,
$0\le k\le n-1$, then obviously $\cd(D_P^k)\subset\ch_k(\R_+,A),
1\le k\le n$, and hence $P$ is $n$--regular.\end{proof}

Of considerable importance are those projections which are $n$--regular
for all $n\in\Z_+$. 

\begin{dfn}\label{I.S1-5.4} Let $P\in\OpAn$ be an orthogonal
projection. The pair $(P,\Ppl)$ is called \emph{elliptic} if 
$P$ is $n$--regular for all $n\in\Z_+$.
% or, equivalently, if \myref{I.G2-5.2} holds for all $s\in\R$.
\end{dfn}
In Proposition \plref{I.S2-5.3} we will express ellipticity completely
in terms of the orthogonal projections $P,\Ppl$ such that it
could equally well be defined for any pair $(P,Q)$ of orthogonal
projections in $\OpAn$.

Now consider an orthogonal projection $P\in\OpAn$ 
such that $(P,\Ppl)$ is elliptic.
Can $k\in \Z_+$ in Proposition \plref{I.S2-4.7} 
be replaced by any nonnegative real number?
We are going to show that if $D_P$ is self--adjoint,
then this is indeed true and follows from 
complex interpolation. 
But since we will need the argument again below
we state the result in the more general framework of
scales of Hilbert spaces:

\begin{dfn}\label{I.S5-4.10}
Let $(H_s)_{s\in\R_{(+)}}$ be a scale of Hilbert spaces 
and let
$B\in\operatorname{Op}^\mu:=\operatorname{Op}^\mu((H_s)_{s\in\R_{(+)}})$.
$B$ is called \emph{regular at $s_0\ge 0$} if
\begin{align}
  &u\in H_0,\; Bu\in H_{s_0} \mbox{ or } B^tu\in H_{s_0} 
     \Longrightarrow u\in H_{{s_0}+\mu}.\label{I.G2-4.14}
\end{align}
If the scale is parametrized over $\R$ and satisfies axiom \myref{Hscalefive}
(cf. Definition \plref{I.S3-2.5}), then 
$B$ is called \emph{elliptic} if for all $s\in\R$ we have
\begin{align}
  &u\in H_{-\infty} ,\; B u \in H_s \mbox{ or } B^t u \in H_s
     \Longrightarrow u\in H_{s+\mu} .
   \label{I.G3-4.13}
\end{align}
\end{dfn}

For the definition of $B^t$ see Proposition \plref{I.S1-2.2} (1).

\begin{prop}\label{I.S2-4.8} Let $(H_s)_{s\in\R}$ be a scale
of Hilbert spaces satisfying axiom \myref{Hscalefive}.
Let $B\in \Op^\mu$ for some  $\mu\ge 1$ and assume that
$B$ is self--adjoint as an operator in $H_0$. 
Then the following statements are equivalent:
\begin{romanenum}
\item
$B$ is regular at all $n\in\Z_+$,
\item
$
 (B-i)^{-1} \in \Op^{-\mu}$,
\item $B$ is elliptic.
\end{romanenum}

If we have only
$\mu>0$ then \textnormal{(ii)} and \textnormal{(iii)} are still
equivalent.
\end{prop}
\begin{remark}\label{I.S5-4.12}\begin{enumerate}
\item If one and hence all of the three equivalent conditions
are fulfilled, then we infer in particular that the domain of $B$,
considered as a self--adjoint operator in $H_0$, is $H_\mu$.
\item Let $B\in\Op^\mu, \mu\ge 0,$ be symmetric, i.e. $\scalar{Bu}{v}=
\scalar{u}{Bv}$ for $u,v\in H_\infty$. If $B$ is regular at $0$,
then $B:H_1\to H_0$ is a self--adjoint operator in $H_0$. Since
$H_\infty$ is dense in $H_1$, $B\restr H_\infty$ is essentially
self--adjoint.
\item If $H_s=H_s(A)$ then the operators $A$ and $|A|^\alpha, \ga\ge 0$,
are elliptic. This follows from the previous proposition and the
fact that for $s\in\R$ the operators
\[\begin{split}
  &(I+|A|^2)^{-s+1/2}(A-i)^{-1}(I+|A|^2)^{s},\\
  &(I+|A|^2)^{-s+\ga/2}(|A|^\ga-i)^{-1}(I+|A|^2)^{s}
  \end{split}
\]
are bounded in $H_0$. The latter follows from the Spectral Theorem.
\end{enumerate}
\end{remark}
\begin{proof}
(iii) $\Rightarrow$ (i) is clear.

(i)$\Rightarrow$ (ii):
Let $v\in H_n , n\in\Z_+.$ Then $u:=(B\pm i)^{-1} v\in H_0$ and thus
\[
  Bu = v\mp i u \in H_0 .
\]
Since $ \mu\ge 1$, the regularity at $0$
implies $u\in H_1$ and iterating this argument shows
$v\pm i u\in H_n$. 
Then (i) implies $u\in H_{n+\mu}.$ 

We have proved that $(B\pm i)^{-1}$ maps $H_n$ into 
$H_{n+\mu}$ for all $n\in\Z_+$. 
Now (ii) follows from Proposition
\plref{I.S1-2.2}. 

(ii)$\Rightarrow$ (iii): We note first that $(B\pm i)^{-1}$ maps
$H_{-\infty}$ bijectively onto $H_{-\infty}$ with inverse $B\pm i$.
Namely, since $B$ is self--adjoint, $(B\pm i)(B\pm i)^{-1}\restr H_0=
\id_{H_0}$ and $(B\pm i)^{-1}(B\pm i)\restr H_1=\id_{H_1}$.
Since $H_\infty\subset H_1\subset H_0$ and $(B\pm i)^{\pm 1}$ leaves
$H_\infty$ invariant, we have $(B\pm i)^{\pm 1}(B\mp i)^{\mp 1}\restr H_\infty
=\id_{H_\infty}$ and by duality these identities also hold for $H_{-\infty}$.

Let $u\in H_{-\infty}$ and $v:=Bu\in H_s.$ 
We have $u\in H_t$ for some $t\in\R$ and
\[
  u=(B-i)^{-1}v-i(B-i)^{-1} u\in H_{\min(s,t)+\mu},
\]
from (ii). 
Iterating this argument we find $u\in H_{s+\mu}.$
\end{proof}

Now assume that $D_P$ is self--adjoint. Then the previous result does
not apply directly to $D_P$ since $D_P$ is not an operator of order
$1$ with respect to the scale $(\ch_s(\R_+,A))_{s\in\R_+}$:
the problem is that
$D_P$ is not defined on $\ch_\infty(\R_+,A)$. However, the proof
of the previous result shows that the following slight modification
of Proposition \plref{I.S2-4.8} is true.

\begin{custthmnum}{\ref{I.S2-4.8}'}
\begin{prop}\label{I.S2-4.8'} Let $(H_s)_{s\in\R_+}$ be a scale
of Hilbert spaces and let $B\in \Op^\mu ,\mu\ge 1$. Let $\tilde B$ be
a self--adjoint operator in $H_0$ with $B\restr\cd(\tilde B)=\tilde B$.
Then the following statements are equivalent:
\begin{romanenum}
\item For all $n\in\Z_+$\textnormal{:}
$u\in\cd(\tilde B), \; Bu\in H_n \Longrightarrow u\in H_{n+\mu}$,
\item
$
 (\tilde B-i)^{-1} \in \Op^{-\mu}$,
\item for all $s\in\R_+$\textnormal{:}
$u\in\cd(\tilde B), \; Bu\in H_s \Longrightarrow u\in H_{s+\mu}.$
\end{romanenum}
\end{prop}
\end{custthmnum}

In view of \myref{I.G2-2.12}, $D$ may be considered as an operator
of order $1$ with respect to the scale of Hilbert spaces
$(\ch_s(\R_+,A))_{s\in\R_+}$. Then Proposition \plref{I.S2-4.8'}
applies to $D_P$.
Summing up we have proved the following regularity theorem.

\begin{theorem}[Regularity Theorem] \label{I.S5-4.13}
Let $P\in\OpAn$ be a $\gamma$--symmetric
orthogonal projection. The pair $(P,\Ppl)$ is elliptic if and only if 
for all $s\in\R_+$ the relations
\begin{equation}\begin{split}
  &D^*f=g\in\ch_s(\R_+,A),\quad f\in L^2(\R_+,H),\\
  &Pf(0)=0,
   \end{split}
   \label{I.G2-4.14a}
\end{equation}
imply $f\in \ch_{s+1}(\R_+,A)$.
\end{theorem}
\begin{proof} The "if" part follows from Proposition \plref{I.S2-4.7}.
For the "only if" part it remains to note that $D_P$ is
self--adjoint, in view of Corollary \plref{I.S3-4.4}. Hence the 
previous discussion applies.\end{proof}

The Regularity Theorem allows to improve Corollary \plref{I.S2-2.14} (2). 

\begin{prop}\label{I.S4-4.9}  For all $s\ge 0$ we have
\begin{equation}\begin{split}
    \ch_{s+1}&(\R_+,A)\\
     &=\bigsetdef{f\in L^2(\R_+,H)}{Df\in \ch_s(\R_+,A),\; f(0)\in H_{s+1/2}}.
      		\end{split}
    \label{I.G2-4.15}
\end{equation}
\end{prop}
\begin{proof}
Let $f$ be in the right hand side of \myref{I.G2-4.15}. Then by
Proposition \plref{I.S1-2.11} we choose $g\in\ch_{s+1}(\R_+,A)$ with 
$g(0)=f(0)$ and put $f_1:=f-g$. Then $Df_1\in\ch_s(\R_+,A), Pf_1(0)=0$
and hence by the Regularity Theorem we have $f_1\in\ch_{s+1}(\R_+,A)$ and
thus $f\in\ch_{s+1}(\R_+,A)$.
\end{proof}

\comment{
This argument relies on the existence of 
a $P$ which is $n$--regular for all $n\in\Z_+$. 
In Proposition \plref{I.S1-5.2} below we will see that 
$\Ppl$ has this property.}%comment

This argument relies on the existence of 
at least one elliptic pair $(P,\Ppl)$.
In Proposition \plref{I.S1-5.2} below we will see that 
the pair $(\Ppl,\Ppl)$ indeed is elliptic.

Let $P\in\OpA{\mu}$ be an orthogonal projection such that
$D_P$ is self--adjoint. Then $\cd(D_P^n)\subset \ch_n(\R_+,A)$ already
implies that $P$ is $n$--regular. Namely, from complex interpolation
we infer for $0\le k\le n$
\begin{equation}\begin{split}
    \cd(D_P^k)&=[\ch_0(\R_+,A),\cd(D_P^n)]_{k/n}\\
       &\subset [\ch_0(\R_+,A), \ch_n(\R_+,A)]_{k/n}\\
      &=  \ch_k(\R_+,A).
  \end{split}
\end{equation}

Corollary \plref{I.S3-4.4} asserts that
$D_P$ is self--adjoint if and only if $P$ is regular and 
$\gamma$--symmetric.
By Theorem \plref{I.S1-4.4} (2), $D_{P,\max}=D_P$ is also
self--adjoint in this case. We now want to show that
$D_{P,\max}$ may be self--adjoint also for non--regular $\gamma$--symmetric
projections. We do not treat the description of self--adjoint extensions
in full generality,
we have singled out only the most tractable
class of projections.

\begin{prop}\label{I.S1-4.7}  Let $P\in\OpAn$ be a $\gamma$--symmetric
orthogonal projection.
If there is $C>0$ such that for $\eps\in(0,1)$
\begin{equation}
\|(I-P) (I+\eps^2 A^2)^{-1} P\|_{H_{-1/2}\to H_{1/2}}\le C,
\label{I.G4-4.17}
\end{equation}
then $D_{P,\max}$ is self--adjoint.
\end{prop}
\begin{proof} In view of Theorem \plref{I.S1-4.4} (1), $D_{P,\max}$ is a closed
operator with $D_{P,\max}\supset \ovl{D_P}=D_{P,\max}^*$. Hence 
we need
to show that $D_{P,\max}$ is symmetric.
 Let $f,g\in \cd(D_{P,\max})$. Then using Lemma \plref{I.S1-4.1}
we find
\[\begin{split}
   -(D_{P,\max}&f,g)+(f,D_{P,\max}g)\\
  &=\lim_{\eps\to 0}
    \bigBss{-1/2}{(I+\eps^2 A^2)^{-1} \gamma  (I-P) f(0)}{%
      (I-P) g(0)}\\
  &=\lim_{\eps\to 0}\bigBss{-1/2}{(I-P) (I+\eps^2 A^2)^{-1} P\gamma f(0)}{%
     g(0)}.
  \end{split}
\]
In $H_0$ the family $(I+\eps^2 A^2)^{-1}$ converges to $I$ strongly, as $\eps \to
0$. Since $(I+\eps^2 A^2)^{-1}$ commutes with $|A|$ we also
have $\lim\limits_{\eps\to 0}(I+\eps^2A^2)^{-1}=I$ strongly in $H_s$, for
each $s\in\R$. In particular, 
for $\xi\in H_{1/2}$ we have
$\lim\limits_{\eps\to 0} (I-P) (I+\eps^2 A^2)^{-1} P\xi=0$  in $H_{1/2}$. Hence
\[ (I-P) (I+\eps^2 A^2)^{-1} P\in\cl(H_{-1/2},H_{1/2})\] converges to $0$
as $\eps\to 0$,
pointwise on a dense subset of $H_{-1/2}$. Consequently, the uniform norm
bound implies the strong convergence, and $D_{P,\max}$ is symmetric.\end{proof}

We note some criteria for $ (I-P) (I+\eps^2 A^2)^{-1} P\in
\cl(H_{-1/2},H_{1/2})$ to be bounded. 
First, boundedness is obvious if
$P$ and $A^2$ commute. More generally, we have

\begin{lemma}\label{I.S4-4.10} If the commutator $[P,A^2]$ is in $\OpA{1}$, then
\myref{I.G4-4.17} holds and hence $D_{P,\max}$ is self--adjoint.
\end{lemma}
\begin{proof} Using the identity
\[(I-P)(I+\eps^2A^2)^{-1}P=
    \eps^2 (I-P) (I+\eps^2A^2)^{-1} [P,A^2](I+\eps^2A^2)^{-1},\]
we estimate
\begin{align}
   &\|(I-P)(I+\eps^2A^2)^{-1}P\|_{H_{-1/2}\to H_{1/2}}\notag\\
   &\le\|I-P\|_{H_{1/2}\to H_{1/2}} \eps^2 \|[P,A^2] (I+\eps^2A^2)^{-1}\|_{H_{-1/2}\to H_{1/2}}\notag\\
   &\le  \|I-P\|_{H_{1/2}\to H_{1/2}} 
  \|[P,A^2]\|_{H_{3/2}\to H_{1/2}}\eps^2\|(I+\eps^2A^2)^{-1}\|_{H_{-1/2}\to
   H_{3/2}}\notag\\
&\le C.\tag*{\qed}
\end{align}
\exendproof

If we restrict to $P\in\cp_s(\PsiA)$, with $\PsiA$ an algebra satisfying
\myref{I.G3-3.22}, then $[P,A^2]\in\OpA{1}$ is implied by the following
condition prominent in Alain Connes' Noncommutative Differential
Geometry:
\begin{MLlist}
\item $(H,|A|,\PsiA)$ forms a spectral triple.\puteqnnum\label{I.G3-4.18}
\end{MLlist}
Indeed, from $[|A|,B]$ bounded
for $B\in\PsiA$, we deduce $[|A|,B]\in\OpAn$ as in \cite[Lemma 1]{Con:GSPV}, 
wherefrom we easily derive
$[P,A^2]\in\OpA{1}$. \myref{I.G3-4.18} will also be important in deriving
the heat asymptotics in part IV of this work. This is plausible
from the fact that for $A$ a classical pseudodifferential operator
with scalar principal symbol on a compact manifold, the algebra
of classical pseudodifferential operators of order $0$ satisfies 
\myref{I.G3-4.18}.

We record the result.

\begin{prop}\label{I.S3-4.7} If $(H,|A|,\PsiA)$ is a spectral triple, then
$D_{P,\max}$ is self--adjoint for any $P\in\cp_s(\PsiA)$.
\end{prop}

Finally, we discuss the existence of non--regular $P\in\OpAn$.

\begin{prop}\label{I.S1-4.8} $D_{\Pmi,\max}$ is self--adjoint; but if $A$ is 
unbounded,
then $\Pmi$ is not regular.
\end{prop}
\begin{proof} Since $[\Pmi,A^2]=0$, the
self--adjointness of $D_{\Pmi,\max}$ follows from Lemma \plref{I.S4-4.10}. 

Since $A$ is unbounded and anticommutes with $\gamma$ the operator
$P_{\ge 1}(A)$ is unbounded, too. Hence 
we can find $u=P_{\ge 1}u\in H\setminus H_{1/2}$. 
Then the function
\[f(x):=e^{-Ax}u\]
is in $\cd(D_{\Pmi,\max})$ but $f(0)\not\in H_{1/2}$; in view
of \myref{I.G3-4.5b}, $\Pmi$ cannot be regular. 
\end{proof}

%Proposition \plref{I.S1-5.10} below contains a more general result.

More generally, we can prove:

\begin{prop}\label{I.S1-5.10} 
Let $A$ be unbounded. 
If $P$ is regular and $[P,A]\in\OpAn$,
then $I-P$ is not regular.
However, in this case $D_{I-P,\max}$ is self--adjoint.
\end{prop}
\begin{proof}
Let $f\in\cd(D^*)$. Since $[P,A]\in\OpAn$
one checks that $Pf, (I-P)f\in\cd(D^*)$. Hence
$Pf\in\cd(D_{I-P,\max}), (I-P)f\in\cd(D_{P,\max})$. If $I-P$ were regular
we would get $\cd(D^*)\subset \ch_1(\R_+,A)$.
In view of the Trace Theorem \plref{I.S1-2.10} this contradicts
\myref{I.G1-4.5}.

The self--adjointness of $D_{I-P,\max}$ follows from Lemma
\plref{I.S4-4.10} since $[P,A^2]=[P,A]A+A[P,A]\in \OpA{1}$.
\end{proof}

\section{Criteria for regularity}
\label{Secfive}

We now study more closely the notion of $n$--regularity established
in Definition \plref{I.S1-4.6}.

\begin{prop}\label{I.S1-5.2} $\Ppl$ is $n$--regular for all $n\in\Z_+$ or,
in other words, the pair $(\Ppl,\Ppl)$ is elliptic.
\end{prop}
\begin{proof} For $n=1$ this follows from Lemma \plref{I.S1-4.3}.

Inductively, we assume that $\Ppl$ is $n$--regular. 
Let $f\in\cd(D_{\Ppl}^{n+1})$, then
$D_{\Ppl}f\in\cd(D_{\Ppl}^n)\subset \ch_n(\R_+,A)$ and consequently
\begin{equation}
     (D_{\Ppl}^jf)(0)\in H_{n+1/2-j},\quad 1\le j\le n+1.
\end{equation}
In view of Corollary \plref{I.S2-2.14}
it remains to prove that $f(0)\in H_{n+1/2}$. 
Since $f\in\cd(D_{\Ppl}^n)\subset \ch_n(\R_+,A)$
we already have $f(0)\in H_{n-1/2}$.
In view of Corollary \plref{I.S1-2.12} we may choose
$h:=e^{(n+1)}(0,(Df)(0),\ldots,(D^nf)(0))
\in\ch_{n+1}(\R,A)$ and put $f_1:=f-h$.

As in \myref{I.G1-4.6} we put
$g(x):=   \chi(x) e^{-Ax}f_1(0), \quad x\le 0$.
Then $g\in\cd((D_-^*)^{n+1})$ and $(D_-^jg)(0)=(D^jf_1)(0), 0\le j \le n$.
Hence Proposition \plref{I.S2-2.13} implies $f_1\in\ch_{n+1}(\R_+,A)$
and thus $f\in\ch_{n+1}(\R_+,A)$.
\end{proof}

Now we are ready to state the analogue of Theorem \plref{I.S1-4.4}
for $n$--regular projections:

\begin{prop}\label{I.S2-5.2} 
Let $P\in\OpAn$ be an orthogonal projection.
The following statements are equivalent:
\begin{romanenum}
\item $P$ is $n$--regular.
\item $r_+^{(k)}(\cd(D_{P,\max}^{k+1}))\subset  
     \moplus\limits_{j=0}^{k}H_{k-j+1/2},\quad 0\le k \le n-1$.
\item For $1\le k\le n$ the following holds:
\begin{MLlist}
\item if $\xi\in H_{-1/2}$ and $P\xi=0, P_- \xi\in H_{k-1/2}$,
    then $\xi\in H_{k-1/2}.$\puteqnnum\label{I.G2-5.1}
\end{MLlist}
\end{romanenum}
\end{prop}

\begin{proof} The equivalence of (i) and (ii) follows from
Corollary \plref{I.S2-2.14}.

(ii)$\Rightarrow$(iii):
We proceed by induction on $n$.
For $n=1$ the assertion follows from Theorem \plref{I.S1-4.4}.

Let $P$ be $n$--regular. It remains to check the case
$k=n$ in \myref{I.G2-5.1}. Consider $\xi\in H_{-1/2}, P\xi=0,
P_-\xi\in H_{n-1/2}.$ The induction hypothesis implies
$\xi\in H_{n-3/2}$. In view of Corollary \plref{I.S1-2.12}
choose $g:=e^{(n-1)}(P_-\xi,0,\ldots,0)\in \ch_n(\R_+,A)$
and put
\[f(x):=\chi(x) e^{-Ax}P_+ \xi+P_- g(x).\]
Then $f\in\cd((D_P^*)^n)$ and hence by the $n$--regularity of $P$, $D_P^*=D_P$ 
and we have $f\in\ch_n(\R_+,A), \xi=f(0)\in H_{n-1/2}$.

(iii)$\Rightarrow$ (i): Again we proceed by induction on $n$.
For $n=1$ the assertion follows from Theorem \plref{I.S1-4.4}.

Assume \myref{I.G2-5.1} for $1\le k\le n$.
By the induction hypothesis $P$ is $(n-1)$--regular, so for
$f\in \cd(D_P^n)$ we have $Df\in\cd(D_P^{n-1})\subset\ch_{n-1}(\R_+,A)$.
Put $f_1:=P_-f$, then $Df_1=P_+Df\in\ch_{n-1}(\R_+,A), P_+f_1(0)=0$.
In view of Proposition \plref{I.S1-5.2} and Proposition
\plref{I.S2-4.7} we find $P_-f(0)=f_1(0)\in H_{n-1/2}$.
Thus \myref{I.G2-5.1} implies $f(0)\in H_{n-1/2}$ and
hence from Corollary \plref{I.S2-2.14} we infer
$f\in\ch_n(\R_+,A)$.\end{proof}

Next we clarify the relation between the notion ``ellipticity''
and ``elliptic pair'' (Definitions \plref{I.S1-5.4} and  \plref{I.S5-4.10}). 

\newcommand{\IAviert}{(I+A^2)^{1/4}}
\newcommand{\IAmviert}{(I+A^2)^{-1/4}}

\begin{prop}\label{I.S2-5.3} Let $P\in \OpAn$ be a $\gamma$--symmetric
orthogonal projection.
\begin{thmenum}
\item $P$ is $n$--regular if and only if the operator
$K:=\IAviert TT^*\IAviert,$ 
where $T= P\Ppl+ (I-P)\Pmi$ (cf. \myref{I.G5-3.5}),
is regular
at $0,1,\ldots, n-1$ in the sense of Definition \plref{I.S5-4.10}.
In particular, the pair $(P,\Ppl)$ is elliptic if and only if
$K$ is elliptic.
\item The pair $(P,\Ppl)$ is elliptic if and only
if for all $s\in\R$ the following holds:
\begin{MLlist}
\item If  $\xi\in H_{-\infty}$ and
 $P\xi=0, \Pmi \xi\in H_s$, then $\xi\in H_s$.
      \puteqnnum\label{I.G2-5.2}
\end{MLlist}
\end{thmenum}
\end{prop}
\begin{proof}
We remark first that \myref{I.G2-5.1} and
\myref{I.G2-5.2} are, in fact, symmetric in $P$ and $P_-$. 
Namely, \myref{I.G2-5.2} is clearly a consequence of
\begin{custeqnnum}{\ref{I.G2-5.2}'}
\begin{MLlist}
\item  $\xi\in H_{-\infty}, P\xi\in H_t, \text{ and } P_-\xi\in H_s 
\text{ imply } \xi\in H_{\min(s,t)}.$\puteqnnum
    \label{I.G2-5.2'}
\end{MLlist}
\end{custeqnnum}
But applying \myref{I.G2-5.2} to $\eta:=\xi-P\xi$ shows that
\myref{I.G2-5.2'} also follows \myref{I.G2-5.2}, so both statements are
equivalent; a similar reasoning works for \myref{I.G2-5.1}.

By the $\gamma$--symmetry of $P$ and $P_+$, \myref{I.G2-5.1} resp. 
\myref{I.G2-5.2} also hold
with $I-P$ and $\Ppl$ in place of $P$ and $P_-$.

With these preparations we prove (1):
Let $K$ be regular at $0,\ldots, n-1$. In view of Proposition
\plref{I.S2-5.2} (iii), we consider $\xi\in H_{-1/2}$ with $P\xi=0$
and $\Pmi\xi\in H_{k-1/2}$ for some $1\le k\le n$. We put
$\eta:=\IAmviert\xi\in H_0$ and find
\[     K\eta=\IAviert (I-P)\Pmi\xi\in H_{k-1}.\]
Since $K$ is regular at $k-1$ we conclude $\eta\in H_k$ and
thus $\xi\in H_{k-1/2}$.

Conversely, let $P$ be $n$--regular and consider $\xi\in H_0$ with
$K\xi\in H_k$ for some $0\le k\le n-1$. $K\xi\in H_k$ implies
in view of Remark \plref{I.S5-4.12} 
\[P\Ppl P \eta, (I-P)\Pmi(I-P)\eta\in H_{k+1/2},
  \quad \eta:=\IAviert\xi.\]
We invoke the $n$--regularity of $P$ and infer from
\[P(\Ppl P\eta)\in H_{k+1/2},\quad \Pmi(\Ppl P\eta)=0\]
that $\Ppl P\eta\in H_{k+1/2}$. Proceeding in this way
we arrive at $P\eta\in H_{k+1/2}$. Analogously one derives
$(I-P)\eta\in H_{k+1/2}$, hence $\eta\in H_{k+1/2}$
and finally $\xi\in H_{k+1}$. Hence $K$ is $n$--regular.

Thus the first assertion of (1) is proved. The second assertion
of (1) as well as (2)
follow from the first one and Proposition \plref{I.S2-4.8}.
\end{proof}

\comment{
To prove (2) we use (1) and
note that $P$ is $n$--regular for all $n\in\Z_+$ if and only if
the self--adjoint operator $H$ is regular at all $n\in\Z_+$.
Now the assertion follows from Proposition \plref{I.S2-4.8}.
\end{proof}}
 
\comment{\begin{dfn}\label{I.S1-5.4} Let $P\in\OpAn$ be a $\gamma$--symmetric orthogonal
projection. $P$ is called {\em elliptic} if 
$P$ is $n$--regular for all $n\in\Z_+$ or, equivalently, if
\myref{I.G2-5.2} holds for all $s\in\R$.
\end{dfn}

In view of Proposition \plref{I.S2-4.8} and Proposition \plref{I.S2-5.3},
$P$ is elliptic if and only if $H=|A|^{1/2}TT^*|A|^{1/2}$ is
elliptic.} % end comment

From now on we will use the axioms \myref{I.G3-3.22} on the algebra
$\PsiA$.

\begin{theorem}\label{I.S1-5.6} Let $P\in\PsiA$ be a $\gamma$--symmetric
orthogonal projection. If the pair $(P,\Ppl)$ is Fredholm, then
it is elliptic.
\end{theorem}
\begin{proof}  Assume first that $(P,\Ppl)$ is an invertible pair.
Then 
\[T=P\Ppl+(I-P)\Pmi\in\PsiA\]
is invertible (by Remark \plref{I.S1-3.9} 
and Corollary \plref{I.S1-3.4}).
By axiom  \myref{I.G3-3.22a} we then have $T^{-1}\in\PsiA$.
Now let $\xi\in H_{-\infty}, \Pmi \xi=0, P \xi\in H_s$. Then
\begin{equation}
   \xi= T^{-1}T\xi = T^{-1}P\Ppl \xi = T^{-1}P\xi\in H_s.
   \label{I.G1-5.8}
\end{equation}
If $P\in\cp(\PsiA)$ is arbitrary, then we write 
$P=P'+P''$ as in  \myref{I.G1-3.10}. The proof of Theorem
\plref{I.S1-3.10} shows that $P'+\Ppl''\in\cp_s(\PsiA)$;
by \myref{I.G1-3.23} the pair $(P'+\Ppl'',\Ppl)$ is invertible.
By the first part of this proof, $(P'+\Ppl'',\Ppl)$ also is
an elliptic pair. By
Lemma \plref{I.S1-2.4} the finite rank operators $P''$ and
$\Ppl''$ 
are smoothing, hence
$(P'+\Ppl'')\xi=P\xi+(\Ppl''-P'')\xi \in H_s$
and as before we see that $\xi\in H_s$.\end{proof}

If $A$ is a discrete operator, one can get a better result. First
we state the abstract elliptic regularity theorem
for scales of Hilbert spaces:

\begin{prop}\label{I.S3-5.7} Let $(H_s)_{s\in\R}$ be a scale of Hilbert
spaces such that for  $s'>s$ the embedding $H_s'\hookrightarrow H_s$
is compact. Let $T\in\Op^\mu:=\Op^\mu((H_s)_{s\in\R}), \mu>0$. Assume that $T$ is regular
at $0$. Then
\begin{thmenum}
\item For $0\le s\le\mu$, $T$ extends to a Fredholm operator
$H_s\to H_{s-\mu}$ with index independent of $s$.
\item Let $S$ be the generalized inverse of $T$ defined by
$ST\xi=\xi$ for $\xi\in\ker T^\perp,$ and $S\xi:=0$ for $\xi\in \im T^\perp$,
where $\perp$ is taken in $H_0$. Then $S$ extends by continuity
to a bounded linear operator $H_s\to H_{s+\mu}, -\mu\le s\le 0$.
\end{thmenum}
If, in addition, $T$ is elliptic then $S$ is a parametrix in the sense that
\[ I-TS, I-ST\in\Op^{-\infty},\]
and for all $s\in\R$,  $T$ extends to a Fredholm operator
$H_s\to H_{s-\mu}$.
\end{prop}
\begin{proof} We first consider $T$ as an unbounded operator in $H_0$;
then $T$ is a closed operator with domain $H_\mu$. To see this
let $(x_n)_{n\in\N}\subset H_\mu$ be a sequence such that
$x_n\to x, Tx_n\to \xi$ in $H_0$. Since $T\in\Op^\mu$ we have
the equality $Tx=\xi$ in $H_{-\mu}$, and the regularity of $T$ at
$0$ implies $x\in H_\mu$. 

By the same argument, $T^t$ induces a closed operator with domain
$H_\mu$, which is the adjoint, $T^*$, of $T:H_\mu\to H_0$.

By assumption, $H_\mu\hookrightarrow H_0$ is compact. Since the
domains of $T$ and $T^*$ are both compactly embedded,
$T$ and $T^*$ are Fredholm operators $H_\mu\to H_0$. By duality,
$T$ and $T^t$ induce Fredholm operators $H_{0}\to H_{-\mu}$ hence,
by complex interpolation, from $H_s\to H_{s-\mu}, 0\le s\le \mu$.

By definition, the generalized inverse $S$ maps $H_0\to H_\mu$ and
its adjoint, $S^*$, is a generalized inverse of $T^*$.
By duality and complex interpolation, again, $S$ and $S^*$ induce
continuous maps $S,S^t:H_s\to H_{s+\mu}, -\mu\le s\le 0$.

The operators $I-ST, I-TS$ are
orthogonal projections in $H_0$ with $(I-ST)(H_0)=\ker T\subset H_\mu,
(I-TS)(H_0)=\ker T^t\subset H_\mu$.
Since $I-TS, I-ST$ are
orthogonal projections, by duality, they map $H_{-\mu}\to H_0$ and hence
$H_{-\mu}\to H_\mu$.
This proves that $\ker T\restr H_s$ is independent of $s$ for 
$0\le s\le \mu $.

If $T$ is elliptic then the finite rank operators
$I-TS$ and $I-ST$ map $H_0\to H_\infty$ and Lemma \plref{I.S1-2.4}
implies $I-TS, I-ST\in\Op^{-\infty}$.
\end{proof}

Now we can state the main result of this section:

\begin{theorem}\label{I.S1-5.8} Let $A$ be discrete. 
For a $\gamma$--symmetric orthogonal projection
$P\in\PsiA$ the following statements are equivalent:
\begin{romanenum}
\item $P$ is regular.
\item $(P,\Ppl)$ is an elliptic pair.
\item $(P,\Ppl)$ is a Fredholm pair.
\end{romanenum}
\end{theorem}
\begin{proof} (i)$\Rightarrow$ (iii): 
Let $P$ be a regular projection. By Proposition \plref{I.S2-5.3}, this
is equivalent to the fact that $K:=\IAviert TT^*\IAviert,$ with 
$T= P\Ppl+ (I-P)\Pmi,$
is regular at $0$. From Proposition \plref{I.S3-5.7} we infer
that $K$ is Fredholm $H_s(A)\to H_{s-1}(A), 0\le s\le 1$. Hence
$TT^*$ is Fredholm $H_s\to H_s, -\frac 12 \le s\le \frac 12$,
in particular, $(P,P_+)$ is a Fredholm pair.
\\
(iii)$\Rightarrow$ (ii): This follows from Theorem \plref{I.S1-5.6}.
\\
(ii)$\Rightarrow$ (i): This follows from Theorem
\plref{I.S1-4.4} and the definition of ellipticity.
\end{proof} 

We note that the implications
"(i)$\Rightarrow$(iii)" and "(ii)$\Rightarrow$(i)" remain valid
if we assume only that $P$ is a $\gamma$--symmetric
orthogonal projection in $\OpAn$.
 
\begin{example} In \cite[Sec. 3]{BruLes:EICNLBVP} we described a special
class of orthogonal projections defining 
generalized Atiyah--Patodi--Singer boundary value problems.
We now discuss these projections in some detail:

Let $P$ be a $\gamma$--symmetric orthogonal projection.
Assume that
\begin{MLlist}
\item $[P,|A|]=0$,\eqnnumlbl{I.G1-5.10}
\item $PAP=\ga |A| P$ for some $\ga>-1$.\eqnnumlbl{I.G1-5.11}
\end{MLlist}
\myref{I.G1-5.10} implies that $P\in\OpAn$.

\begin{prop}\label{I.S1-5.9}
Let $P$ be a $\gamma$--symmetric orthogonal projection
satisfying \myref{I.G1-5.10} and \myref{I.G1-5.11}.
Then the following assertions hold:
\begin{thmenum}
\item The pair $(P,\Ppl)$ is elliptic.
\item If $\dim \ker A<\infty$ then $(P,\Ppl)$ is
a Fredholm pair. If $\ker A=0$ then $(P,\Ppl)$ is an
invertible pair.
\end{thmenum}
\end{prop}
\begin{proof} For $\gl>0$ consider the operator
\begin{equation}
        P_\gl:= \frac 12 (|A|+\gl)^{-1}(A+|A|).
\end{equation}
$P_\gl$ converges to $P_{>0}(A)$ strongly, as $\gl \to 0+$.
In view of \myref{I.G1-5.10}, \myref{I.G1-5.11} we have
\begin{equation}
   \begin{split}
      PP_\gl P&= (|A|+\gl)^{-1}|A| \frac{1+\ga}{2} P\\
              &\to (P_{<0}(A)+P_{>0}(A))\frac{1+\ga}{2}P
   \end{split}
\end{equation}
strongly, as $\gl\to 0+$. Hence
\begin{align}
    PP_{>0}(A)P&=(P_{<0}(A)+P_{>0}(A))\frac{1+\ga}{2}P,\\
\intertext{and consequently we infer that}
    P\Ppl P&= \frac{1+\ga}{2}P+R,
   \label{I.G3-5.8}\\
\intertext{where}
    R&= P(\Ppl-P_{>0}(A))P- \frac{1+\ga}{2}P_0P\in\Op^{-\infty}(A).
\end{align}
By construction, 
if $\dim\ker A<\infty$,
then $R$ is of finite rank and if $\ker A=0$ then $R=0$.
For $T:= PP_++(I-P)P_-$ we find
\begin{equation}
    TT^*=\frac{1+\ga}{2} I +R+\gamma R \gamma^*.
\end{equation}
Since $R$ is smoothing, $TT^*$ and thus the pair $(P,\Ppl)$ is elliptic,
in view of Proposition
\plref{I.S2-5.3}.
If  $\dim\ker A<\infty$ then $TT^*$ and thus $(P,P_+)$ 
is Fredholm. Finally, if  $\ker A=0$ then $TT^*$
and hence $(P,P_+)$ is invertible.\end{proof}
\end{example}
\begin{prop}\label{I.S1-5.11} Let $P\in\OpAn$ be a $\gamma$--symmetric 
orthogonal projection.
\begin{thmenum}
\item If $(I-P)\Ppl(I-P)\in \OpA{{-\eps}}$ for some $\eps>0$, then
the pair $(P,\Ppl)$ is elliptic.
\item Let $P_1\in\OpAn$ be another $\gamma$--symmetric orthogonal projection.
If $(P,\Ppl)$ is elliptic and $(P-P_1)\in\OpA{-\eps}$ for some 
$\eps>0$, then $(P_1,\Ppl)$ is also elliptic.
\end{thmenum}
\end{prop}
\begin{proof} (1)
Let $x\in H_{s_0}$, for some $s_0\in \R, Px=0, \Pmi x\in H_s$. We have
\[0=Px=P\Ppl x+P\Pmi x,\]
thus $P\Ppl x=-P\Pmi x\in H_s$. This implies
\begin{align*}
    \Ppl x&=P\Ppl x+(I-P)\Ppl(I-P)x\\
      &=-P\Pmi x+(I-P)\Ppl(I-P)x\\
      &\in H_{\min(s,s_0+\eps)}.
\end{align*}
Iterating this procedure we conclude $x\in H_s$.

(2) Let $x\in H_{s_0}, P_1x=0, \Pmi x\in H_s$. Then
\[Px=(P-P_1)x\in H_{s_0+\eps},\]
hence $x\in H_{\min(s,s_0+\eps)}$.
Again, we we reach the conclusion by iterating this
argument.\end{proof}

Now we can present the

\begin{proof}[Proof of Theorem \plref{intro-S2}]  In this proof $D_P$ denotes
the operator defined in \myref{I.G1-4.2}, i.e.
$\cd(D_P)=\bigsetdef{u\in \ch_1(\R_+,A)}{Pu(0)=0}$,
and we put
$D_{P,0}:=D_P\restr\cd_P:=\bigsetdef{f\in C^\infty_0 (\R_+, H_1)}{P f(0)=0}$.
Note that the $D_P$ in the formulation of Theorem \plref{intro-S2} 
is the $D_{P,0}$ in the present notation.

By Corollary \plref{I.S1-4.5} $\cd_P$ is a core for 
$D_P$ and hence $D_P$ is the closure of $D_{P,0}$.
Furthermore, in view of \myref{intro-1.13} and 
Theorem \plref{I.S1-4.4} (1) we have $D_P^*=D_{P,\max}.$
From \myref{intro-1.20} and Theorem \plref{I.S1-5.6} we infer
that $P$ is a regular projection and hence in view of the
three equivalent characterizations of regularity in
Theorem \plref{I.S1-4.4} (2) we have $D_P=D_{P,\max}=D_P^*$.
This proves the first part of Theorem \plref{intro-S2}.

To prove the second part we assume that $A$ is discrete
and that $D_P$ is self--adjoint on
$\bigsetdef{f\in\ch_1(\R_+,A)}{Pf(0)=0}$. 
Hence $P$ is regular by Theorem \plref{I.S1-4.4}.
Then \myref{intro-1.20}
is a consequence of Theorem \plref{I.S1-5.8}. \myref{intro-1.13}
follows from Corollary \plref{I.S3-4.4}.
\end{proof}

\section{Variable coefficients and second order operators}
\label{Secsix}

We now study the model operator \myref{intro-1.10} with variable coefficients. 
In most applications we will need the
following results only in a small neighborhood of $x=0$. Therefore, we
may assume that the operator has constant coefficients at infinity.
More precisely, we consider the differential operator
\begin{equation}\begin{split}
  D&=\gamma(\frac{d}{dx}+A(x))\\
   &=:\gamma(\frac{d}{dx}+A_0)+x\gamma A_1(x)\\
   &=:D_0+x\gamma A_1(x),
		\end{split}
  \label{I.G1-6.1}
\end{equation}
where \myref{intro-1.9} now holds with $A(x)$ in place of $A$, for all $x$, 
and we assume in addition that
\begin{alphnumbering}[I.G5-6.2]
%\begin{MLnumlist}
%\newlength{\boxwidth}
%\setlength{\boxwidth}{\textwidth}
%\addtolength{\boxwidth}{-\leftmargin}
%\addtolength{\boxwidth}{-1.5cm}
\begin{MLlist}
\begin{MLnumlist}
\MLitem{$A(x)\in \Op^1(A_0)$ is self--adjoint with domain
$H_1(A_0)$ and elliptic with respect to the scale $(H_s(A_0))_{s\in\R}$,}
\eqnnumlbl{I.G5-6.2a}

\medskip
%\begin{itemize}%\renewcommand{\labelitemii}{$\bullet$}
\MLitem{for all $s\in\R$, the map}
\begin{equation}
\R\ni x\mapsto (I+A_0^2)^s A(x)(I+A_0^2)^{-1/2-s}\in\cl(H_0(A_0))
   \label{I.G5-6.2b}
\end{equation}
\item is smooth,%}\eqnnumlbl{I.G5-6.2b}

\medskip
\item   $A_1(x)=0$ for $|x|\ge R>0$.\puteqnnum\label{I.G5-6.2c}
\end{MLnumlist}
\end{MLlist}
Note that, from \myref{intro-1.9}
\begin{MLlist}
\item  $\gamma A_1(x)+A_1(x)\gamma=0$, \quad for all $x$.
\puteqnnum\label{I.G5-6.2d}
%\end{itemize}
\end{MLlist}
\end{alphnumbering}
\myref{I.G5-6.2a} implies that
\begin{equation}
H_s(A(x))=H_s(A_0)\quad\text{for } x,s\in\R,
\end{equation}
and \myref{I.G5-6.2b} 
implies
\begin{equation}
A^{(j)}(x)\in\Op^1(A_0),\quad\text{for } j\in\Z_+.
\end{equation}
We abbreviate $H_s:=H_s(A_0)$.
As in the constant coefficient case we are interested in $D$ as an 
unbounded operator in $L^2(\R_+,H)$ with domain
$\cinfz{(0,\infty),H_\infty}$; 
the corresponding operator in $L^2(\R,H)$ with domain 
$\cinfz{\R,H_\infty}$ will be denoted by $\tilde D$.
Then $D, \tilde D$ are symmetric operators in the respective Hilbert
spaces. We have
\begin{equation}
     D^2=-\frac{d^2}{dx^2}+A(x)^2-A'(x).
    \label{I.G1-6.4}						      
\end{equation}
Our analysis of the constant coefficient case 
applies to $D_0$ in \myref{I.G1-6.1}. 
In particular, $\tilde D_0$ is essentially
self--adjoint on $\cinfz{\R,H_\infty}\subset L^2(\R,H)$.
Its closure has domain $\ch_1(\R,A_0)$ and will also be denoted by 
$\tilde D_0$. Moreover, $\ch_1(\R,A_0)\subset\cd(\tilde D)$.

We first state the analogue of Lemma \plref{I.S1-4.1}.

\begin{lemma}\label{I.S1-6.1} 
\begin{thmenum}
\item $\cd(D^*)\subset \{f\in L^2(\R_+,H)\,|\, f'\in L^2(\R_+,H_{-1})\}$
and we have a continuous restriction map
$r:\cd(D^*)\to H_{-1/2}, f\mapsto f(0)$.
\item For $f,g\in \cd(D^*)$ we have
  \[-(D^*f,g)+(f,D^*g)=
  \lim_{\eps\to 0}\Big(\gamma (\pm\eps A_0+i)^{-1} f(0), 
    (\eps A_0+i)^{-1} g(0)\Big).\]
If $f(0)$ or $g(0)$ lies in $H_{1/2}$ then
\begin{equation}
 -(D^*f,g)+(f,D^*g)=\Bss{\mp 1/2}{\gamma f(0)}{g(0)}.
  \label{I.G1-6.5}
\end{equation}
\end{thmenum}
\end{lemma}
\begin{proof} The proof is almost identical to the proof of Lemma
\plref{I.S1-4.1}. We only have to note that
for $f\in \cd(D^*)$ the family
\[f_\eps(x):=i (\eps A(x)+i)^{-1}f(x)\]
satisfies the relations $f_\eps\in L^2(\R_+,H)$,
$\lim\limits_{\eps\to 0} f_\eps =f$,
\[ f_\eps'(x)= - \eps (\eps A(x)+i)^{-1} A'(x)f_\eps +
         i (\eps A(x)+i)^{-1} f'(x)\in L^2(\R_+,H),\]
and $D^*f_\eps$ has uniformly bounded norm.
Thus $D^*f_\eps\rightharpoonup D^*f$ as shown by the next lemma.\end{proof}

The last proof used the following simple but useful lemma, which allows 
to describe a core of an unbounded operator in terms of weak convergence.

\begin{lemma}\label{I.S1-6.2} Let $D$ be a closed operator in a Hilbert space
$\ch$.
\begin{thmenum}
\item Let $f\in\ch$ and $(f_n)\subset \cd(D)$ be a sequence 
%with
%$\lim\limits_{n\to\infty}\|f-f_n\|_\ch=0$ 
such that $(f_n)$ converges weakly to $f\in \ch$ and $(Df_n)$ 
has uniformly bounded norm in $\ch$.
Then $f\in \cd(D)$ and $Df_n\rightharpoonup Df$.
\item Let $\ce\subset\cd(D)$ be a linear subspace such that for
every $f\in \cd(D)$ there is a sequence $(f_n)\subset \ce$ with
$f_n\rightharpoonup f$ and $(Df_n)$ bounded. Then $\ce$ is a core for $D$.
\end{thmenum}
\end{lemma}
\begin{proof} (1) For every $g\in\cd(D^*)$ we have
\def\MLto{\mathop{\longrightarrow}}
\begin{equation}
   \label{I.G1-6.6}
   (Df_n,g)=(f_n,D^*g)\MLto_{n\to\infty} (f,D^*g),
\end{equation}
i.e. the bounded sequence $(Df_n)$ converges weakly on a dense subspace
and hence is weakly convergent. Let $h$ be the weak limit of $(Df_n)$.
Then in view of \myref{I.G1-6.6} we have $(h,g)=(f,D^*g)$ for all $g\in\cd(D^*)$,
thus $f\in\cd(D)$ and $Df=h$.

(2) To prove that $\ce$ is a core for $D$ we have to show that
if $g\in\cd(D)$ with $0=(g,f)+(Dg,Df)$ for all $f\in \ce$ then $g=0$.
Given such a $g$ we choose $(f_n)\subset \ce$ with $f_n\rightharpoonup g$ and
$Df_n\rightharpoonup Dg$, using (1). Then
\begin{equation}
(g,g)+(Dg,Dg)=\lim_{n\to\infty}\{ (g,f_n)+(Dg,Df_n)\}=0.\tag*{\qed}
\end{equation}
\exendproof

As in the constant coefficient case we want to find boundary
conditions defining self--adjoint extensions of $D$. Our main tool
will be (a variant of) the Kato--Rellich Theorem so we need to estimate 
relative bounds.

\begin{lemma} \label{I.S5-6.4}
Let $B(x), x\in\R$, be a family of operators of order
$1$ satisfying \myref{I.G5-6.2b}. Then for $\eps>0$ and $n\in\Z_+$ there
exists $C(\eps,n,B)$ such that for $f\in\ch_{n+1}(\R,A_0)$ 
we have the estimate
\[  \begin{split}\bigl\|& Bf\bigr\|_{\ch_n}\\
    & \le
     c_n\Bigl( \sup_{x\in\R}\max\bigl(\bigl\| B(x) T^{-1}\bigr\|_{\cl(H_0)},
      \bigl\| T^nB(x) T^{-1-n}\bigr\|_{\cl(H_0)}\bigr)
      +\eps\Bigr) \|f\|_{\ch_{n+1}}\\
    &\quad+ C(\eps,n,B) \|f\|_{\ch_0}.
    \end{split}
\]
Here, $c_n$ is a universal constant depending only on $n$ and
$T:=(I+A_0^2)^{1/2}$.  
\end{lemma}
\begin{proof} In this proof,  $c_n$ and $C(\cdot)$ denote generic
constants depending on their respective arguments.

For $f\in\ch_1(\R,A_0)$ we have
\begin{equation}\begin{split}
     \|Bf\|_{\ch_0}^2&=\int_\R \|B(x) f(x)\|_{H_0}^2 dx\\
      &\le \sup_{x\in\R}\| B(x)T^{-1}\|_{\cl(H_0)}^2 \;
        \int_\R \|Tf(x) \|_{H_0}^2 dx\\
   &= \sup_{x\in\R}\| B(x)T^{-1}\|_{\cl(H_0)}^2 \;
            \|Tf\|_{\ch_0(\R,A_0)}^2.
\end{split}\label{I.G5-6.8}
\end{equation}
Using this estimate and Proposition \plref{I.S5-2.10} we find
\begin{align*}
\bigl\| &B f\bigr\|_{\ch_n}\\
     &\le c_n \sum_{j=0}^n \bigl\|\bigl(\frac{d}{d x}\bigr)^j T^{n-j} Bf   \bigr\|_{\ch_0}\\
     &\le c_n \sum_{j=0}^n \sum_{i=0}^j 
          \bigl\|T^{n-j} B^{(i)}f^{(j-i)}\bigr\|_{\ch_0}\\
     &\le c_n \sum_{j=0}^n \sum_{i=0}^j 
          \sup_{x\in\R}\bigl\| T^{n-j}B^{(i)}(x)T^{j-n-1}\bigr\|_{\cl(H_0)}
          \bigl\|T^{n+1-j} f^{(j-i)}\bigr\|_{\ch_0},
\end{align*}
where in the last inequality we have used
\myref{I.G5-6.8} with $T^{n-j}B^{(i)}T^{j-n}$ in place of $B$ and
$T^{n-j}f^{(j-i)}$ in place of $f$. We treat the summands with
$i=0$ and $i\not=0$ separately. Summands with $i=0$ are estimated
using the following inequality which follows from complex interpolation:
\begin{alphnumbering}
\begin{align}
\|T^{n-j}&B(x)T^{j-n-1}\|_{\cl(H_0)}\notag\\
   &\le \max\bigl( 
        \|B(x)T^{-1}\|_{\cl(H_0)}, \|T^n B(x)T^{-n-1}\|_{\cl(H_0)}\bigr).
      \label{I.G5-6.9a}
\intertext{For $i>0$ we estimate in view of Proposition \plref{I.S5-2.10}
and \myref{I.G2-2.12}}
   \bigl\|T^{n+1-j} f^{(j-i)}\bigr\|_{\ch_0}&\le 
            c_n\bigl\|f^{(j-i)}\bigr\|_{\ch_{n+1-j}}\nonumber\\
      &\le c_n \bigl\|f\bigr\|_{\ch_{n+1-i}} \nonumber\\
      &\le \eps \|f\|_{\ch_{n+1}}+ C(\eps,n)\|f\|_{\ch_0}.
\end{align}
The last inequality follows from $i>0$, the fact that 
$\ch_s(\R,A)=H_s(\tilde D)$ \myref{I.G1-2.10}, and
the Spectral Theorem. The lemma is proved.\end{alphnumbering}\end{proof}

Next we note a variant of the Kato--Rellich Theorem.

\begin{prop}\label{I.S5-6.5} Let $S$ be a self--adjoint operator
in the Hilbert space $H$ and let $R\in\Op^1(S)$ be a symmetric
operator. Assume that for fixed $N\in\Z_+$ and $\xi \in H_{\infty}(S)$
we have the estimates
\begin{alphnumbering}[I.G5-6.10]
\begin{align}
     \|R\xi\|_{H_0(S)}&\le b \|S\xi\|_{H_0(S)}+c\|\xi\|_{H_0(S)}
        \label{I.G5-6.10a}, \\
     \|R\xi\|_{H_{N-1}(S)}&\le b \|S\xi\|_{H_{N-1}(S)}+c\|\xi\|_{H_0(S)},
        \label{I.G5-6.10b}
\end{align}
with $b<1$ and $c>0$.
\end{alphnumbering}

Then the operator $S+R$ is self--adjoint with domain $H_1(S)$
and regular at $0\le k\le N-1$.
\end{prop}
\begin{proof} The estimates \myref{I.G5-6.10} show that for $0<\eps<1-b$
there is a $\gl_0>0$ such that for
$\gl\ge \gl_0$ 
\begin{equation}
    \bigl\|R (S-i\gl)^{-1}\|_{\cl(H_s(S))}\le b+\eps<1.
     \label{I.G5-6.11}
\end{equation}
A priori, \myref{I.G5-6.11} holds for $s=0$ and $s=N-1$, but
by complex interpolation it holds for $0\le s\le N-1$.
Thus the resolvent
\begin{equation}
 (S+R-i\gl)^{-1}= (S-i\gl)^{-1} 
\sum_{n=0}^\infty (-1)^n \bigl(R(S-i\gl)^{-1}\bigr)^n
\label{I.G5-6.12}
\end{equation}
converges in $\cl(H_s(S),H_{s+1}(S)), 0\le s\le N-1$.
As in the proof of the implication "(i)$\Rightarrow$(ii)" of
Proposition \plref{I.S2-4.8} we now conclude that $S+R$ 
is regular at $0\le s\le N-1$.
Of course, the self--adjointness of $S+R$ with domain $H_1(S)$
also follows from \myref{I.G5-6.12} (as in the proof
of the Kato--Rellich Theorem).
\end{proof}

Now we can prove the main result of this section,
the generalization of the Regularity Theorem \plref{I.S5-4.13}
to variable coefficients.
 
\begin{theorem}[Regularity Theorem]\label{I.S5-6.6}
Assume \myref{I.G1-6.1}--\myref{I.G5-6.2}.
\begin{thmenum}
\item The operator $\tilde D$ is elliptic with respect to the
scale of Hilbert spaces $(\ch_s(\R,A_0))_{s\in\R}$. In particular,
$\tilde D$ is essentially self--adjoint and $\cd(\ovl{\tilde D})=
\ch_1(\R,A_0).$

\item 
Let $P\in\OpAn$ be a $\gamma$--symmetric
orthogonal projection. 
\begin{enumerate}
\item \newcommand{\locali}{\textnormal{(i)}} The pair $(P,\Ppl)$ is elliptic if and only if
for all $s\in\R_+$ the relations
\begin{equation}\begin{split}
  &D^*f=g\in\ch_s(\R_+,A),\quad f\in L^2(\R_+,H),\\
  &Pf(0)=0,
   \end{split}
   \label{I.G5-6.13}
\end{equation}
imply $f\in \ch_{s+1}(\R_+,A)$. 

In other words, in this case the operator $D_P$ defined by
\[D_P:= D^*\restr \cd(D_{0,P})=\{f\in \ch_1(\R_+,A_0)\,|\, Pf(0)=0\}\]
fulfills the equivalent conditions \textup{(i)--(iii)} of
Proposition \plref{I.S2-4.8'}.
In particular, $D_P$ is self--adjoint.
\item $P$ is regular if and only if the operator $D_P$ defined
in \locali\ is self--adjoint. 
\end{enumerate}
\end{thmenum}
\end{theorem}
\begin{proof} 
(1) By Proposition \plref{I.S2-4.8} and Remark
\plref{I.S5-4.12} 2., it suffices to prove the regularity 
at all $n\in\Z_+$. 
Fix $x_0\in\R$ and
write
\[A(x)=A(x_0)+(x-x_0)A_1(x_0,x).\]
From \myref{I.G5-6.2b} we infer $A_1(x_0,x)\in\Op^1(A_0)$ and
in view of \myref{I.G5-6.2a} we have
$H_s(A(x_0))=H_s(A_0)$ and $\Op^\mu(A(x_0))=\Op^\mu(A_0)$
for all $s,\mu\in\R$.

In view of Lemma \plref{I.S5-6.4} we may
choose
$\varphi,\psi\in\cinfz{\R}$ with $\varphi\equiv 1$ in a neighborhood
of $x_0$ and $\psi\equiv 1$ in a neighborhood of $\supp\varphi$ such
that for the operator 
\[B:=(x-x_0)\psi(x)\gamma A_1(x_0,x)\]
we have, for a fixed $n$,
\begin{equation}\begin{split}
    \norm{Bf}_{\ch_{n}}&\le \tilde b \norm{f}_{\ch_{n+1}}+C(n)\norm{f}_0\\
      &\le  \tilde b c_n \norm{\gamma(\frac{d}{dx}+A(x_0))f}_{\ch_{n}}
                +C'(n)\norm{f}_0,
    		\end{split}
     \label{I.G6-6.14}
\end{equation}
and
\begin{equation}
    \norm{Bf}_0\le \tilde b \norm{f}_{\ch_1}+C(0)\norm{f}_0\le
        \tilde b c_0 \norm{\gamma(\frac{d}{dx}+A(x_0))f}_0+C'(0)
                    \norm{f}_0, 
     \label{I.G6-6.15}
\end{equation}
with $\tilde b$ so small that $b:=\tilde b \max( c_0,c_n) <1$.
Hence $B(x)$ is
$\gamma(\frac{d}{dx}+A(x_0))$--bounded with relative bound
$b<1$ for both the $\ch_0$ and the $\ch_n$--norm.

Now 
consider $f\in \ch_0(\R,A_0)$ with $\tilde D f=g\in \ch_n(\R,A_0)$. 
By Lemma \plref{I.S1-2.8} and Proposition \plref{I.S5-6.5}
the operator
\begin{equation}
     \tilde D(x_0):=\gamma(\frac{d}{dx}+A(x_0)+(x-x_0)\psi(x)A_1(x_0,x))
        \label{I.G1-6.10}
\end{equation}
is self--adjoint with domain $H_1(A_0)$
 and regular at all $0\le s\le n$.

Assume that we already know that $f\in\ch_t(\R,A_0)$.
Then
\begin{equation} 
   \tilde D(x_0)^*(\varphi f)=\tilde D^*(\varphi f)=\gamma \varphi' f+
         \varphi \tilde D^*f\in \ch_{\min(n,t)},
    \label{I.G1-6.11}
\end{equation}
and consequently $\varphi f\in \ch_{\min(n,t)+1}(\R,A_0)$.

Since $x_0$ was arbitrary we have proved
$f\in \ch_{\min(n,t)+1,\loc}(\R,A_0)$.

To prove that $f\in \ch_{\min(n,t)+1}(\R,A_0)$ we choose a cut--off function
$\varphi\in\cinfz{\R}$ with $\varphi\restr[-R,R]\equiv 1$. Then
$\varphi f\in\ch_{\min(n,t)+1}(\R,A_0)$ and in view of 
\myref{I.G5-6.2c} we have
$(1-\varphi)f\in\cd(\tilde D_0^{n+1})=\ch_{n+1}(\R,A_0)$.
(1) is proved.

(2i) We assume that $(P,\Ppl)$ is elliptic and fix an integer $n\ge 0$. 
By the Regularity Theorem \plref{I.S5-4.13},
the operator $D_{0,P}$, with $D_0$ from \myref{I.G1-6.1},
induces a closed operator in $\ch_n(\R_+,A_0)$
with domain $\cd(D_{0,P})\cap \ch_{n+1}(\R_+,A_0)$. Hence, for
$f\in\cd(D_{0,P})\cap \ch_{n+1}(\R_+,A_0)$ there is an estimate

\begin{equation}\begin{split}
     \norm{f}_{\ch_{n+1}(\R_+,A_0)}
         &\le c\bigl( \norm{D_{0,P} f}_{\ch_{n}(\R_+,A_0)}+ 
              \norm{f}_{\ch_{n}(\R_+,A_0)}\bigr)\\
         &\le c  \norm{D_{0,P} f}_{\ch_{n}(\R_+,A_0)}+\frac 12 
                \norm{f}_{\ch_{n+1}(\R_+,A_0)}+c_{1/2}\norm{f}_0,
		\end{split}
     \label{I.G7-6.18}
\end{equation}
thus
\begin{equation}
      \norm{f}_{\ch_{n+1}(\R_+,A_0)}
            \le c\bigl(\norm{D_{0,P} f}_{\ch_{n}(\R_+,A_0)}+\norm{f}_0\bigr).
      \label{I.G7-6.19}
\end{equation}
In view of Lemma \plref{I.S5-6.4} we can thus choose 
cut--off functions $\varphi, \psi$ near $0$ as in the first part of this
proof (cf. \myref{I.G6-6.14}, \myref{I.G6-6.15})  such that
$x \psi(x) \gamma A_1(x)$ is $D_{0,P}$--bounded with relative bound $b<1$,
simultaneously for the $\ch_0(\R_+,A_0)$ and the
$\ch_{n}(\R_+,A_0)$--norm.

Proposition \plref{I.S5-6.5} does not apply directly with $R=x\psi(x)\gamma
A_1(x)$ and $S=D_{0,P}$ since we do not necessarily have
$R\in\Op^1(D_{0,P})$. However, an inspection of the proof of Proposition
\plref{I.S5-6.5} shows that
\[  \tilde D_{0,P}:=D_{0,P}+x\gamma \psi(x) A_1(x)\]
is self--adjoint with domain $\ch_1(\R_+,A_0)$ and for $0\le t\le n$
the following holds
(cf. Proposition  \plref{I.S2-4.8'} versus Proposition \plref{I.S2-4.8}):
\[ f\in\cd(\tilde D_{0,P}^*), \; \tilde D_{0,P}^*f\in \ch_t(\R_+,A_0)
   \Longrightarrow f\in \ch_{t+1}(\R_+,A_0).\]
Now consider $f$ satisfying \myref{I.G5-6.13}.
Assume we already know $f\in\ch_t(\R_+,A_0)$. As in \myref{I.G1-6.11}
we conclude $\varphi f\in\ch_{\min(n,t)+1}(\R_+,A_0)$. Using part
(1) of the theorem we infer from
\[ \tilde D^*(1-\varphi) f= -\gamma \varphi' f+(1-\varphi) D^* f\in
     \ch_{\min(n,t)}(\R,A_0)\]
that also $(1-\varphi)f\in\ch_{\min(n,t)+1}(\R_+,A_0)$.

Conversely, assume that \myref{I.G5-6.13} implies $f\in \ch_{s+1}(\R_+,A_0)$.
Then we reverse the roles of $D_P$ and $D_{0,P}$; i.e. \myref{I.G5-6.13}
implies that $D_P$ induces a closed operator in $\ch_n(\R_+,A_0)$
with domain $\cd(D_{0,P})\cap \ch_{n+1}(\R_+,A_0)$. Hence, the estimates
\myref{I.G7-6.18} and \myref{I.G7-6.19} hold with $D_P$ in place 
of $D_{0,P}$. Then we choose cut--off functions $\varphi,\psi$ as
before such that $x\gamma \psi(x) A_1(x)$ is $D_P$--bounded with
relative bound $b<1$, simultaneously for the $\ch_0(\R_+,A_0)$ and
the $\ch_n(\R_+,A_0)$--norm; as before $n\ge 0$ is some fixed integer.

Then we see again that the operator
\[ D_{P,\psi}:= D_P - x\gamma \psi(x) A_1(x)= 
     D_{0,P}+x(1-\psi(x))\gamma A_1(x)\]
satisfies
\begin{equation}
 f\in\cd(\tilde D_{P,\psi}^*), \; \tilde D_{P,\psi}^*f\in \ch_t(\R_+,A_0)
   \Longrightarrow f\in \ch_{t+1}(\R_+,A_0)
   \label{I.G7-6.20}
\end{equation}
for $0\le t\le n$. 
Next consider $f\in L^2(\R_+,H)$ with
\[ D_0^*f=g\in \ch_s(\R_+,A_0), \quad Pf(0)=0.\]
From part (1) we infer $f\in\ch_{s+1,\loc}((0,\infty),A_0)$
and from 
\begin{equation}\begin{split}
       D_{P,\psi}^*\varphi f &= \gamma \varphi' f+ \varphi D_{P,\psi}^* f\\
            &= \gamma \varphi' f + \varphi D_0^* f\in \ch_s(\R_+,A_0),
\end{split}\label{I.G7-6.21}\end{equation}
we infer, in view of \myref{I.G7-6.20}, that $\varphi f\in
 \ch_{\min(n,s)+1}(\R_+,A_0)$ and hence $f\in \ch_{\min(n,s)+1}(\R_+,A_0)$.

Since $n$ is arbitrary we have proved that the operator $D_{0,P}$
satisfies \myref{I.G2-4.14a} and hence from Theorem \plref{I.S5-4.13}
we conclude that the pair $(P,\Ppl)$ is elliptic.

Specializing this proof to $s=0$ immediately implies (2ii).
The theorem is proved.
\end{proof}

\begin{remark} 
An alternative proof of the first part of this theorem could
have been given by means of an operator valued pseudodifferential
calculus as in \cite[Sec. 2]{BruSee:RESORSO}.
\end{remark}

The regularity at $0$ or $1$ of $\tilde D$ resp. $D_P$ can be obtained
under weaker assumptions. We record the results and leave the details
to the reader.

\begin{theorem}\label{I.S5-6.7} If we replace \textup{(\ref{I.G5-6.2a},b)}
by
\begin{MLlist}
\item $A(x)\in\Op^1(A_0)$ and $x\mapsto A(x)(I+A_0^2)^{-1/2}$ is
continuous,
\end{MLlist}
then the following hold.
\begin{thmenum}
\item
$\tilde D$ is essentially self--adjoint and $\cd(\ovl{\tilde D})=
\ch_1(\R,A_0);$ in particular, $\tilde D$ is regular at $0$.
\item
Let $P\in\OpAnn$ be a regular $\gamma$--symmetric orthogonal
projection. Then $D$ is self--adjoint on
\[\cd(D_P):=\cd(D_{0,P})=\{f\in \ch_1(\R_+,A_0)\,|\, Pf(0)=0\}.\]
\end{thmenum}
\end{theorem}

\begin{theorem}\label{I.S5-6.8} If we replace \textup{(\ref{I.G5-6.2a},b)}
by
\begin{MLlist}
\item $A(x), A'(x)\in\Op^1(A_0)$ and $x\mapsto A(x)(I+A_0^2)^{-1/2}$ is
in $C^1(\R,\cl(H_0))$,
\end{MLlist}
then the following hold.
\begin{thmenum}
\item $\tilde D^2$ is essentially self--adjoint and $\cd(\ovl{\tilde D^2})=
\ch_2(\R,A_0).$ In particular, $\tilde D$ is regular at $0$ and $1$.
\item
Let $P\in\OpAnn$ be a $2$--regular $\gamma$--symmetric orthogonal
projection. Then 
\[\cd(D_P^2)=\cd(D_{0,P}^2)=\{f\in \ch_2(\R_+,A_0)\,|\, Pf(0)=0,
\; P(Df)(0)=0\}.\]
\end{thmenum}
\end{theorem}

\section{Well--posed boundary value problems}
\label{Secseven}

In this section we consider the situation described in Sec.\ref{SecOneB}. 
We assume that $A$ is a symmetric
elliptic differential operator of first order
acting on $\cinf{E_N}$, and that $\gamma$ is a bundle endomorphism
such that \myref{intro-1.9} holds. Then
$D$ in \myref{intro-1.10} is an elliptic differential operator on the cylinder $M=\R_+\times N$.
Also, the axioms \myref{I.G3-3.22} are satisfied
with $\Psi^0(A)=\Psi^0_{\text{cl}}(E\restr N).$
For $\xi\in T_x^*(N)$ denote by $N_\pm(\xi)$ the space spanned
by eigenvectors of $\hat A(\xi)$ with positive and negative
eigenvalues, respectively, where $\hat A(\xi)$ denotes the
leading symbol.

According to Seeley \cite[Def. VI.3]{See:TPO},
a classical pseudodifferential operator $P$ of order
$0$ on $\cinf{E_N}$ is called \emph{well--posed} if
\begin{romanenum}
\item $P:H_s(E_N)\to H_s(E_N)$ has closed range for each $s\in\R$;
\item for each $\xi\in T_x^*N\setminus\{0\}$
the principal symbol $\hat P(\xi)$
maps $N_+(\xi)$ injectively onto the range of $\hat P(\xi)$.
\end{romanenum}

$P$ can always be replaced by an orthogonal projection
with the same null space \cite[Lemma VI.3]{See:TPO}, 
hence defining the same boundary condition.
Our aim is to give a functional analytic characterization of well--posedness
for orthogonal projections.

\begin{prop}\label{I.S1-7.1} Let 
$P$ and $Q$ be orthogonal projections in $L^2(E_N)$ which are
classical pseudodifferential operators of order $0$
on $\cinf{E_N}$, and denote
by $\hat P, \hat Q$ their principal symbols. Then $(P,Q)$ is a Fredholm pair
if and only if for each $\xi\in T_x^*M\setminus\{0\}$
\begin{equation}
\hat Q(\xi):\Im \hat P(\xi)\longrightarrow \Im \hat Q(\xi)
\label{I.G1-7.1}
\end{equation}
is an isomorphism.
\end{prop}
\begin{proof} Assume that $(P,Q)$ is Fredholm. Then,
by Prop. \plref{I.S1-3.3}, 
$P-Q\pm 1$ is Fredholm and hence also elliptic
(cf. eg. \cite[Ch. 19.5]{Hor:ALPDOI-IV}).
Thus for $\xi\in T^*N\setminus
\{0\}$ the endomorphisms
\begin{equation}
\hat P(\xi)-\hat Q(\xi)\pm 1
\label{I.G1-7.2}
\end{equation}
are invertible. Now we apply Corollary \plref{I.S1-3.4} to see
that the map \myref{I.G1-7.1} is invertible.

Conversely, if \myref{I.G1-7.1} is invertible then $(\hat P(\xi),
\hat Q(\xi))$ is a finite--dimensional Fredholm pair and
the maps in \myref{I.G1-7.2} are invertible in view
of \myref{I.G1-3.3},
hence $P-Q\pm 1$ is elliptic. Invoking again 
Prop. \plref{I.S1-3.2}, we conclude that $(P,Q)$ is Fredholm.\end{proof}

Condition (i) above is automatic for a pseudodifferential idempotent.
Hence, applying this proposition to the pair $(P,\Ppl)$ immediately
gives

\begin{theorem}\label{I.S1-7.2} Let $P$ be a classical
pseudodifferential operator
of order $0$ which is an orthogonal projection. Then $P$ is well--posed
if and only if $(P,\Ppl)$ is a Fredholm pair.
\end{theorem}

Finally, we give the

\begin{proof}[Proof of Theorem \plref{intro-S4}]
We use the notation of Section \plref{SecOneD}. 
By \myref{intro-1.23}--\myref{intro-1.25}  and Lemma \plref{intro-S1}
we have
\[ \Phi D \Phi^*= \gamma\left(\frac{d}{dx} +A(x)\right)+V(x).\]
We choose a function $\psi\in\cinfz{-\eps_0,\eps_0}, 0\le \psi\le 1,$ with
$\psi=1$ near $0$ such that the operator
\begin{equation}
  A_\psi(x)=\psi(x) A(x) +(1-\psi(x)) A_0
\end{equation}
is elliptic for all $x\ge 0$, which is certainly the case if the
support of $\psi$ is small enough. Then we introduce the
operator
\begin{equation}
   D_\psi:=\gamma\left(\frac{d}{dx} +A_\psi(x)\right)
\end{equation}
and note that $D_\psi$ satisfies \myref{I.G1-6.1} and \myref{I.G5-6.2}.
We will prove the following 

\newtheorem*{claim}{Claim}
\begin{claim} The operator $D_P$ is self--adjoint if and only
if the operator $D_{P,\psi}:=D_\psi^*\restr 
\bigsetdef{f\in\ch_1(\R_+,A_0)}{Pf(0)=0}$ is self--adjoint.
\end{claim}

Theorem \plref{intro-S4} is an immediate consequence of the Claim:

Since the algebra $\Psi^0_{\class} (E_N)$ satisfies
the assumptions \myref{axiomsintro} (cf. the discussion before Theorem
\plref{intro-S2}) and since the elliptic operator $A_0$ is discrete
we infer from Theorem \plref{I.S1-5.8} that the Fredholmness
of the pair $(P,\Ppl)$ is equivalent to the ellipticity of the pair $(P,\Ppl)$
resp. to the regularity of $P$.

From the Regularity Theorem \plref{I.S5-6.6} (2ii) we thus know that
$D_{P,\psi}$ is self--adjoint if and only if the pair $(P,\Ppl)$ is
Fredholm. In view of Theorem \plref{I.S1-7.2} this, in turn,
is equivalent to the fact that $P$ is well--posed in the sense
of Seeley. Again, from the Regularity Theorem \plref{I.S5-6.6} and 
Theorem \plref{intro-S3} we see that the operators 
$Q, K_r, K_l$ are obtained by patching together an interior
parametrix of $D$ and the analogue of \myref{intro-1.21} for
$D_{P,\psi}$.

It remains to prove the claim:

First, we note that the operator $V(x)$ is bounded and hence
adding or subtracting $V(x)$ does not affect the domain of an operator.

Choose a cut--off function $\chi\in\cinfz{-\eps_0,\eps_0}$ such that
$\psi=1$ in a neighborhood of $\supp \chi$. Assume first that
$D_P$ is self--adjoint and let $f\in\cd(D_{P,\psi}^*)$. By the same
calculcation as in \myref{I.G7-6.21} one shows  
$\chi f\in\cd(D_{P,\psi}^*)$, hence by construction we have
$\Phi^*(\chi f)\in \cd(D_P^*)=\cd(D_P)$. Thus, by \myref{intro-1.25},
we have $\chi f\in\cd(D_{P,\psi})$. The proof of the converse is
similar.
The Claim and hence the Theorem are proved.
\end{proof}

%
%  Bibliography stuff
%  the bibliography command unfortunately has some deficiencies
%  - not everybody likes uppercase running heads
%  - it's not included in the table of contents, at least not
%    in article documentclass
%  For use with other documentclasses, you may comment out
%  the following lines
%\addcontentsline{toc}{section}{References}
%\renewcommand\MakeUppercase{\relax}
%\renewcommand\refname{{\upshape References}}
%\bibliography{localbib,mlabbr,localbibs9798,books,papers}%,mathlib}
%\bibliographystyle{lesch}
% Copy of the last bibtex run

 \newcommand{\nop}[1]{} \newcommand{\single}[1]{#1}
  \newcommand{\SwapArgs}[2]{#2#1}
  \newcommand{\translationof}{{English Translation of }}
  \newcommand{\submitted}{{Submitted }}
  \newcommand{\submittedto}{{Submitted to  }} 
  \newcommand{\privcomm}{{Private communication}}
  \newcommand{\toappear}{{to appear}}

\normalfont
Manuscripts can be retrieved from the home page of the second named author.
The papers having a \texttt{math.DG} or \texttt{dg-ga} 
number may be retrieved from \texttt{http://xxx.lanl.gov}.
\signature % redefined in arttitle.tex, for use with article doc class

\end{document}